\documentclass[10pt]{article}
\usepackage[nohead,margin=1.0in]{geometry}
\usepackage{amssymb, amsmath, amsthm,amsfonts}
\usepackage{graphicx,epsfig}
\usepackage[tight]{subfigure}
\usepackage{cite}
\usepackage{times}
\usepackage{algpseudocode}
\usepackage{float}
\usepackage{appendix}
\usepackage{color}
\usepackage{epstopdf}
\usepackage{mathrsfs}
\usepackage{algorithm}
\usepackage{booktabs}
\usepackage{verbatim}
\usepackage[colorlinks,linktocpage,linkcolor=blue]{hyperref}

\graphicspath{{./pics/}}
\usepackage{booktabs}
\usepackage{threeparttable}

\numberwithin{equation}{section}
\newtheorem{theorem}{Theorem}[section]
\newtheorem{remark}{Remark}[section]
\newtheorem{proposition}{Proposition}[section]
\newtheorem{lemma}{Lemma}[section]
\newtheorem{corollary}{Corollary}[section]

\newtheorem{assumption}[theorem]{Assumption}

\allowdisplaybreaks

\def\al{\alpha}

\def\dx{\mathrm{d} x}
\def\dt{\mathrm{d} t}
\def\ds{\mathrm{d} s}
\def\dz{\mathrm{d} z}
\def\dy{\mathrm{d} y}
\def\II{(\Omega)}
\def\la{\lambda}
\def\IPP{IPP}
\def\P01{P_\mathcal{A}}

\title{Numerical Recovery of a Time-Dependent Potential in Subdiffusion\thanks{The work of B. Jin is supported by a start-up fund and Direct Grant of Research, both from The Chinese University of Hong Kong, and Hong Kong General Research Fund (Project No. 14306423). The work of K. Shin is supported by Basic Science Research Program through the National Research Foundation of Korea
(NRF) funded by the Ministry of Education (Grant No. 2019R1A6A1A11051177). The work of  Z. Zhou is supported by Hong Kong Research Grants Council (No. 15303021) and an internal grant of Hong Kong Polytechnic University (Project ID: P0038888, Work Programme: 1-ZVX3).}}

\author{Bangti Jin\thanks{Department of Mathematics, The Chinese University of Hong Kong, Shatin, New Territories, Hong Kong, P.R. China (\texttt{bangti.jin@gmail.com, b.jin@cuhk.edu.hk}).}
\and Kwancheol Shin\thanks{Department of Mathematics, Ewha Womans University, 52, Ewhayeodae-gil, Seodaemun-gu, Seoul
03760, Republic of Korea (\texttt{kcshin3623@gmail.com})}
\and Zhi Zhou\thanks{Department of Applied Mathematics,
The Hong Kong Polytechnic University, Kowloon, Hong Kong, P.R. China (\texttt{ zhizhou@polyu.edu.hk})}\footnotemark[3]}

\begin{document}

\maketitle

\begin{abstract}
In this work we investigate an inverse problem of recovering a time-dependent potential in a semilinear subdiffusion model from an integral measurement of the solution over the domain. The model involves the Djrbashian--Caputo fractional derivative in time. Theoretically, we prove a novel conditional Lipschitz stability result, and numerically, we develop an easy-to-implement fixed point iteration for recovering the unknown coefficient. In addition, we establish rigorous error bounds on the discrete approximation. These results are obtained by crucially using smoothing properties of the solution operators and suitable choice of a weighted $L^p(0,T)$ norm. The efficiency and accuracy of the scheme are showcased on several numerical experiments in one- and two-dimensions.\vskip5pt
\noindent\textbf{Keywords}: inverse potential problem, subdiffusion, Lipschitz stability, error estimate, fixed point method
\end{abstract}

\section{Introduction}
This work is concerned with the recovery of a time-dependent potential coefficient $q(t)$ in a subdiffusion model.
Let $\Omega\subset\mathbb{R}^d $ ($d=1,2,3$) be a convex polyhedral domain with a boundary $\partial\Omega$. Consider the following initial-boundary value problem for the  semilinear  subdiffusion model:
 \begin{equation}\label{eqn:fde}
 \left\{\begin{aligned}
     \partial_t^\alpha u(x,t) -  \nabla \cdot (a(x) \nabla u(x,t)) + q(t)u(x,t) &=f(u(x,t),x,t), &&\mbox{in } \Omega\times{(0,T]},\\
      a(x)\partial_\nu u(x,t)&=g(x,t),&&\mbox{on } \partial\Omega\times{(0,T]},\\
    u(x,0)&=u_0(x),&&\mbox{in }\Omega,
  \end{aligned}\right.
 \end{equation}
where $T>0$ is a fixed terminal time, $f$ and $u_0 \in L^2(\Omega)$ are the given nonlinear source term and initial data, respectively, $\nu$ is the unit outward normal vector to $\partial \Omega$, and the diffusion coefficient $a\in C^2(\overline\Omega)$ satisfies {$a\geq a_0$ for some $a_0>0$}. The fractional order $\alpha \in (0,1)$ is fixed, and the notation $\partial^{\alpha}_t u $ denotes the so-called Djrbashian-Caputo fractional derivative of order $\alpha$ in time, which is defined by \cite{KilbasSrivastavaTrujillo:2006, Jin:2021}
\begin{equation}
    \partial^{\alpha}_t u(t) = \frac{1}{\Gamma (1 - \alpha)} \int_0^t (t-s)^{-\alpha} u'(s){\rm d}s,
\end{equation}
where $\Gamma(z) = \int_0^{\infty} s^{z-1}e^{-s}\ds$ for $\Re (z)>0$ denotes Euler's Gamma function ($\Re$ denotes taking the real part of $z\in \mathbb{C}$). Throughout we denote the solution $u$ corresponding to the potential $q$ by $u(x,t;q)$ or $u(q)$.

The model \eqref{eqn:fde} has received increasing attention over the past three decades, due to the extraordinary capability of the model for describing anomalously slow diffusion processes, also known as subdiffusion. At a microscopical level, it is often described by continuous time random walk, where the waiting time between consecutive jumps follows a heavy-tailed distribution with a divergent mean, in a manner similar to Brownian motion for the standard diffusion equation ($\alpha=1$).
The mathematical model has found many important applications in physics, biology and finance, e.g., thermal diffusion in fractal media \cite{Nigmatulin:1986}, subsurface flow \cite{AdamsGelhar:1992}, and protein transport in membrane \cite{Kou:2008}.
We refer interested readers to the comprehensive review \cite{MetzlerKlafter:2000} for physical motivation of the model and an extensive list of applications in physics and engineering.

In this work, we are concerned with the following {inverse potential problem (IPP)}: given an integral measurement of the solution $u$:
\begin{equation}\label{eqn:m}
    m(t) = \int_{\Omega} u(x,t;q^\dag) \dx,  \quad t \in [0, T],
\end{equation}
we aim to recover the unknown time-dependent potential {$q^\dag(t), \, t\in [0,T]$}. {Integral observations of the type \eqref{eqn:m} have been employed in \cite{PrilepkoOrlovskii:1985,PrilepkoOrlovskii:1985b,CannonLin:1988,CannonLin:1988b,Kamynin:2011,Kamyninukharova:2019,Kamynin:2020} for inverse problems for standard parabolic problems (including degenerate ones).}
Physically, if $u$ denotes the concentration of a chemical diffusing in the domain $\Omega$, then the quantity $\int_\Omega u(x,t) \,\dx$ in \eqref{eqn:m} denotes the total mass of the chemical at time $t$.
Throughout we assume that the ground-truth potential $q^\dag$ belongs to the admissible set ${\mathcal{Q} }= \{ q\in C[0,T]:~ 0\le q \le c_0 \}$ for some $c_0>0$.
In practice, we are interested in the numerical reconstruction of the potential $q(t)$ from the noisy measurement:
\begin{equation*}
    m_\delta(t) =  \int_{\Omega} u(x,t;q^\dag) \dx  + \xi(t), \quad t\in (0,T],
\end{equation*}
and $\xi\in C[0,T]$ denotes the (deterministic) pointwise measurement noise.
The accuracy of the observational data $m_\delta$ is measured by the noise level $\delta$, defined by
$$\|m_\delta - m\|_{C[0,T]} = \delta.$$

In this work we make the following contributions. First, we derive a novel (conditional) Lipschitz stability estimate under minor conditions in Theorem \ref{thm:stab}: for any $p\in [1,\infty]$,
\begin{equation*}
    \|q_1-q_2\|_{L^p(0,T)} \leq c \|\partial_t^\alpha (m_1-m_2)\|_{L^p(0,T)}.
\end{equation*}
This estimate is obtained using smoothing properties of the solution operators in Lemma \ref{lem:op}, and a suitable choice of the weighted $L^p(0,T)$ norm. Second, in Theorem \ref{thm:recon}, we develop an easy-to-implement fully discrete fixed point scheme, and establish the following error bound in the $\ell^p$ norm, $1<p<\infty$, for the discrete approximation $(q_*^n)_{n=1}^N$:
\begin{equation*}
  \|(q_*^n-q^\dag(t_n))_{n=1}^N\|_{\ell^p(\mathbb{R})} \leq c (\delta\tau^{-\alpha}+h^2 + \tau^{1/p} |\log \tau|).
\end{equation*}
This bound is explicit in terms of the discretization parameters $h$ and $\tau$ and the noise level $\delta$ etc, which provides useful guidelines for choosing the algorithmic parameters (i.e., discretization parameters $h$ and $\tau$). This result is achieved using the weighted $\ell^p$ norm as well as suitable error estimates for nonsmooth data. Third and last, we present numerical experiments with both smooth and nonsmooth potentials to illustrate the efficiency of the proposed scheme.

Inverse problems for time-fractional PDEs have received a lot of attention in the last fifteen years and many interesting theoretical and numerical results have been obtained \cite{JinRundell:2015,LiYamamoto:2019}. However, the relevant inverse problems involving a time-dependent elliptic operator are far less studied, since several powerful tools for time-independent elliptic operators, e.g., Laplace transform and Mittag-Leffler functions, cannot be applied directly. Thus, there are only a few results on related linear and nonlinear inverse problems \cite{Zhang:2016,FujishiroKian:2016,SunZhangWei:2019,WeiLiao:2022,WeiXiao:2022,YanZhangWei:2022,YanWei:2023}. Zhang \cite{Zhang:2016} proved the unique recovery of a time-dependent diffusion coefficient $a(t)$ in the model $\partial_t^\alpha u = a(t) u_{xx}$ from the lateral flux observation $-a(t)u_x(0,t)$ using a fixed point argument, and presented numerical results illustrating the approach. Fujishiro and Kian \cite[Theorem 2.2]{FujishiroKian:2016} proved a unique recovery of the time-dependent factor $f(t)$ in the potential term under a suitable positivity condition. {Zhang and Zhou \cite{ZhangZhou:2023} studied the backward problem of the subdiffusion model with time-dependent coefficients, and proved stability for both small and large terminal time scenarios}. In a series of works \cite{SunZhangWei:2019,WeiLiao:2022,WeiXiao:2022,YanZhangWei:2022,YanWei:2023}, Wei et al proved the uniqueness of the recovery of a time-dependent coefficient in the subdiffusion ($\alpha\in(0,1)$) or diffusion-wave ($\alpha\in(1,2)$) model from different types of observations, and presented extensive numerical illustrations of the feasibility of recovery. For example, under a suitable positivity condition, \cite[Theorem 1.1]{SunZhangWei:2019} extended the stability estimate in \cite[Theorem 2.2]{FujishiroKian:2016} to multi-term subdiffusion, which involves multiple fractional orders.
See also \cite{MaSun:2023, Hendy:2022} for the uniqueness results of the related inverse potential problem and inverse source problem, as well as their numerical treatment. However, none of these interesting existing works provides error estimates for the reconstruction algorithms. One goal of the work is to fill this gap of the existing literature.

{Inverse potential problems for the classical parabolic equation have been extensively studied in the literature. Isakov \cite{Isakov:1992} proposed a convergent iterative scheme for recovering a space-dependent potential $q(x)$ from the terminal data under suitable monotonicity conditions on the problem data. Choulli and Yamamoto \cite{ChoulliYamamoto:1996} proved a generic well-posedness result for the inverse potential problem of recovering a space-dependent potential $q(x)$ from the terminal observation in H\"{o}lder spaces. Later, Choulli and Yamamoto \cite{ChoulliYamamoto:1997} proved a Lipschitz stability result for the inverse problem in a Hilbert space setting using the inverse function theorem. This result was shown when the initial data $u_0$ belongs to a suitable subset, and the terminal time $T$ is sufficiently small. These interesting works have been extended and refined in various ways, including the time-fractional case \cite{ZhangZhou:2017,KaltenbacherRundell:2019,JinZhou:2021ip,ZhangZhangZhou:2022,JinKianZhou:2023}. Nonetheless, all these works focus on recovering a space-dependent potential $q(x)$ from terminal observation. The present work follows this long line of research but investigates recovering a time-dependent potential from the integral observation $m(t)$ for a semilinear time-fractional diffusion model. Due to the alignment of the directions of the unknown $q(t)$ with the observational data $m(t)$, the conditional Lipschitz stability estimate does hold, in a manner similar to the case of recovering space-dependent potential $q(x)$ from the terminal data. However, the analysis requires different techniques due to limited smoothing properties of the solution operators, and further, the work provides a complete error analysis of the fully discrete scheme.}

The rest of the paper is organized as follows. In Section \ref{sec:prelim}, we collect preliminary results about the direct problem, e.g., well-posedness and regularity. In Section \ref{sec:stability}, we prove the first main result, i.e., conditional Lipschitz stability of IPP. In Section \ref{sec:recon}, we develop a fully
discrete iterative scheme, and provide a thorough error analysis of the  scheme. Finally, in Section \ref{sec:numer}, we present numerical results to illustrate the performance of the numerical scheme. Throughout, the notation $(\cdot,\cdot)$ denotes the standard $L^2(\Omega)$ inner product and $(\cdot,\cdot)_{L^2(\partial\Omega)}$ the $L^2(\partial\Omega)$ inner product. We often write a bivariate function $f(x,t)$ as $f(t)$ as a vector valued function. The notation $c$ denotes a generic constant which may change at each occurrence, but it is always independent of the noise level $\delta$ and the discretization parameters $h$ and $\tau$, {time step $n$ and iteration index $k$ etc.}

\section{Preliminaries}
\label{sec:prelim}
First we present preliminary results about the solution theory of the model \eqref{eqn:fde}, following the recent textbooks \cite{KubicaYamamoto:2020,Jin:2021}. These results will play a crucial role in establishing the stability and conducting the numerical analysis of the algorithm below. The model involves a nonlinear source and a nonzero Neumann boundary condition, and thus requires slight reworking of the well-posedness result as well as the regularity estimate.
Let $A$ be the $L^2(\Omega)$
realization of the elliptic operator $\mathcal{A} v (x)= -\nabla\cdot(a(x)\nabla v (x))$, with a domain $D(A):=\{v\in
L^2(\Omega): {\mathcal{A} v}\in L^2(\Omega), \partial_\nu v|_{\partial\Omega}=0\}$.
Let $\{\la_\ell\}_{\ell=1}^\infty$ and $\{\varphi_\ell\}_{\ell=1}^\infty$ be, respectively, eigenvalues (ordered nondecreasingly with multiplicity
counted) and the $L^2(\Omega)$-orthonormal eigenfunctions of $A$.  Note that $\la_1=0$ (and has multiplicity $1$)
and the corresponding eigenfunction $\varphi_1=|\Omega|^{-1/2}$
is constant valued, where $|\Omega|$ denotes the Lebesgue measure of the set $\Omega$.

For $u_0 \in L^2\II$ and $f\in L^\infty(0,T;L^2\II)$,
the solution $u$ to the abstract linear evolution problem
 \begin{equation}\label{eqn:fde-0}
 \begin{aligned}
     \partial_t^\alpha u +Au &=f~~\text{for}~ t\in(0,T],\quad \text{with}~~u(0)= u_0,
  \end{aligned}
 \end{equation}
can be represented by
\begin{equation}\label{eqn:sol}
u(t) = F(t)u_0 + \int_0^t E(t-s)f(s)\,\ds,
\end{equation}
where the solution operators $F(t)$ and $E(t)$ are respectively defined by \cite[Section 6.2]{Jin:2021}
\begin{align*}
F(t):=\frac{1}{2\pi {\rm i}}\int_{\Gamma_{\theta,\delta }}e^{zt} z^{\alpha-1} (z^\alpha+A)^{-1}\, \dz \quad\mbox{and}\quad
E(t):=\frac{1}{2\pi {\rm i}}\int_{\Gamma_{\theta,\delta}}e^{zt}  (z^\alpha +A)^{-1}\, \dz ,
\end{align*}
with the integration over a contour $\Gamma_{\theta,\delta}$ in the complex plane $\mathbb{C}$
(oriented counterclockwise), defined by
\begin{equation*}
  \Gamma_{\theta,\delta}=\left\{z\in \mathbb{C}: |z|=\delta, |\arg z|\le \theta\right\}\cup
  \{z\in \mathbb{C}: z=\rho e^{\pm\mathrm{i}\theta}, \rho\ge \delta\} .
\end{equation*}
Throughout, we fix $\theta \in(\frac{\pi}{2},\pi)$ so that $z^{\al} \in \Sigma_{\al\theta}
\subset \Sigma_{\theta}:=\{0\neq z\in\mathbb{C}: {\rm arg}(z)\leq\theta\},$ for all $z\in\Sigma_{\theta}$.
Below we use extensively the following resolvent estimate for the operator $A$:
\begin{equation} \label{eqn:resol}
  \| (z +A)^{-1} \|_{L^2\II \rightarrow L^2\II}\le c_\phi |z|^{-1},  \quad \forall z \in \Sigma_{\phi},
  \,\,\,\forall\,\phi\in(0,\pi).
\end{equation}

The next lemma gives smoothing properties of  $F(t)$ and $E(t)$.
\begin{lemma}\label{lem:op}
For the solution operators $F(t)$ and $E(t)$ defined in \eqref{eqn:sol}, there is some constant $c$ independent of $t$ such that the following estimates hold for all $t>0$.
\begin{itemize}
\item[$\rm(i)$] $\|A F (t)v\|_{L^2\II} +
{t^{1-\alpha}  \| AE (t)v  \|_{L^2\II}} \le c  t^{-\alpha} \|v\|_{L^2\II}$;
\item[$\rm(ii)$] $\|F(t)v\|_{L^2\II} +  t^{1-\alpha}\|E(t)v\|_{L^2\II} \le c   \|v\|_{L^2\II}$.
\end{itemize}
\end{lemma}
\begin{proof}
The estimates in the Dirichlet case can be found in \cite[Theorem 6.4]{Jin:2021}, and the Neumann case follows identically from the resolvent estimate \eqref{eqn:resol}.
\end{proof}

Using these technical tools, we can prove the following existence, uniqueness and regularity of a solution to the nonlinear problem \eqref{eqn:fde}. Due to the nonzero Neumann boundary condition, and the nonlinear source, the result does not appear known. We provide a proof of the result in the appendix for the convenience of readers.
\begin{theorem}\label{thm:sol-reg}
Let $q\in  \mathcal{Q}$ and be a piecewise $C^1$ function, $u_0 \in H^2\II$, $g\in C^1([0,T];H^{\frac12}(\partial\Omega))$, with the compatibility condition $a \partial_\nu u_0 = g(0)$ on $\partial\Omega$, and let $f(u,x,t)$ be smooth and globally Lipschitz continuous.
Then problem \eqref{eqn:fde} admits a unique solution  $u\in C^\alpha([0,T];L^2(\Omega)) \cap
C([0,T];H^2(\Omega))$ and $\partial_t^\alpha u\in C([0,T];L^2(\Omega))$. Moreover,
\begin{align*}
&u\in C^\alpha([0,T];L^2(\Omega)) \cap
C([0,T]; {H^2(\Omega)}) ,\quad
\partial_t^\alpha u\in C([0,T];L^2(\Omega)).
\end{align*}
Moreover, if $q \in \mathcal{Q} \cap C^1[0,T]$, then
\begin{align*}
\partial_tu(t)\in L^2(\Omega)\quad
\mbox{and}\quad
 \|\partial_tu(t)\|_{L^2(\Omega)}\le ct^{\alpha-1}
\quad \mbox{for}\,\,\, t\in(0,T],
\end{align*}
{where the constant $c$ depends on the fractional order $\alpha$, $\| q \|_{C^1[0,T]}$, $\|u_0\|_{H^2(\Omega)}$, $\| g \|_{C^{1}([0,T];H^\frac{1}{2}\II)}$,  Lipschitz constant of $f$ and the terminal time $T$.}
\end{theorem}

\section{Stability and iterative algorithm}\label{sec:stability}
In this section, we establish a novel conditional Lipschitz stability result, and propose a convergent fixed point algorithm. Throughout we make the following assumptions. Condition (i) plays a role in the unique determination and can be ensured by suitable maximum principle (see, e.g., \cite[Section 6.5]{Jin:2021} and \cite{LuchkoYamamoto:2019}), and condition (ii) imposes mild regularity assumptions on the exact data, which holds true under the assumptions in Theorem \ref{thm:sol-reg}.
\begin{assumption}\label{ass:stab}
Let the assumptions in Theorem \ref{thm:sol-reg} hold. Moreover, the following conditions hold.
\begin{itemize}
\item[{\rm(i)}] $m(t)$ is strictly positive and uniformly bounded such that $0 < m_*\le m(t)\le m^*$, {$0\le t\leq T$};
\item[{\rm(ii)}] $m\in C^\alpha[0,T] \cap C^2(0,T]$ such that $|\partial_t^\ell (\partial_t^\alpha m)(t)| \le c t^{-\ell}$ for $\ell=0,1$ and $|\partial_t^{s} m(t)| \le c t^{\alpha-s}$ for $s=1,2$.
\end{itemize}
\end{assumption}

Next, we derive a conditional Lipschitz stability of \IPP. The overall proof strategy is also useful in the analysis of the reconstruction algorithm. Next, for any $\lambda>0$, we define a weighted norm in $L^p(0,T)$ by
\begin{equation}\label{equ: network realization}
 \|  v  \|_{L_\lambda^p(0,T)}^p = \left\{\begin{aligned}
		& \Big(\int_0^T |e^{-\lambda t} v(t)|^p \,\dt\Big)^{\frac1p},&& p\in[0,\infty),\\
		&  \text{esssup}_{t\in(0,T)} e^{-\lambda t}v(t) ,&& p = \infty.
	\end{aligned}\right.
\end{equation}
The $L_\lambda^p(0,T)$ norm is equivalent to the standard $L^p(0,T)$ norm for any $p\in{[1,\infty]}$.
Similarly, for a Banach space $X$, $\| \cdot \|_{L_\lambda^p(0,T;X)}$ denotes the weighted norm of the Bochner space  $L^p(0,T;B)$.
\begin{theorem}\label{thm:stab}
Let $q_i \in {\mathcal{Q}}$ be the solution to \IPP{} with exact data $m_i$, $i=1,2$.
Moreover, let $q_1$ and $m_1$  satisfy Assumption \ref{ass:stab}, and $m_2$ satisfy Assumption \ref{ass:stab} (i).
Then the following stability estimate holds for all $p\in[1,\infty]$,
$$ \| q_1 - q_2  \|_{L^p(0,T)} \le c \| \partial_t^\alpha(m_1  - m_2)   \|_{L^p(0,T)},$$
{where the constant $c$ depends on the bounds $m_*$ and $m^*$  from Assumption \ref{ass:stab}(i), the fractional order $\alpha$,  $\|\partial_t^\alpha m_1(t)\|_{C[0,T]}$, $\|g\|_{C([0,T];L^2(\partial\Omega))}$, and Lipschitz constant of $f$.}
\end{theorem}
\begin{proof}
By integrating over the domain $\Omega$ on both sides of the first equation in \eqref{eqn:fde} and then applying integration by parts in space, we get
\begin{equation}
    \partial^{\alpha}_t m(t) - \int_{\partial\Omega} g(x,t){\rm d}s + q(t) m(t) = \int_{\Omega}f(u, x,t){\rm d}x.
\end{equation}
Note that $m_i(t) \ge m_* >0$ for all $t \in [0,T]$ (Assumption \ref{ass:stab} (i)). Then the potential $q_i$ can be expressed as
\begin{equation*}
    q_i(t) = \frac{\int_{\Omega}f(u(t;q_i),t) \dx - \partial^{\alpha}_t m_i(t) +\int_{\partial \Omega} g(t) \ds  }{m_i(t)} .
\end{equation*}
Consequently, we have
\begin{align*}
 q_1 - q_2  =&  \frac{\int_{\Omega}f(u(t;q_1),t) \dx - \partial^{\alpha}_t m_1(t) +\int_{\partial \Omega} g(t) \ds  }{m_1(t)}
 -\frac{\int_{\Omega}f(u(t;q_2),t) \dx - \partial^{\alpha}_t m_2(t) +\int_{\partial \Omega} g(t) \ds  }{m_2(t)} \\
 =&  \left(\frac{\int_{\Omega}f(u(t;q_1),t) \dx}{m_1(t)}-\frac{\int_{\Omega}f(u(t;q_2),t) \dx}{m_2(t)}\right) + \left({\frac{\partial^{\alpha}_t m_2(t)   }{m_2(t)}} - \frac{\partial^{\alpha}_t m_1(t)}{m_1(t)}\right)\\
 &+ \left(\frac{1  }{m_1(t)} - \frac{1}{m_2(t)} \right)\int_{\partial \Omega} g(t) \ds := \sum_{i=1}^3 \mathrm{I}_i.
\end{align*}
Now we bound the three terms ${\rm I}_i$ using the weighted norm $\|\cdot\|_{L^p_\lambda(0,T)}$. By Assumption \ref{ass:stab}, i.e.,  $0< m_*\le m_i\le m^*$ and $\| u(q_1)  \|_{C([0,T];L^2\II)} \le c$ (in view of Theorem \ref{thm:sol-reg}), and Lipschitz continuity of $f$ in $u$, we derive
\begin{equation*}
\begin{split}
\|   \mathrm{I}_1 \|_{L_\lambda^p(0,T)} & \le \Big\| \frac{ m_1  \int_{\Omega}f(u(t;q_1),t) - f(u(t;q_2),t)\dx }{m_1m_2} \Big\|_{L_\lambda^p(0,T)}
+  \Big\| \frac{ (m_1 - m_2) \int_{\Omega}f(u(t;q_1),t)\dx }{m_1m_2} \Big\|_{L_\lambda^p(0,T)} \\
&\le c \| u(q_1) - u(q_2) \|_{L_\lambda^p(0,T;L^2\II)} + c \| m_1 - m_2 \|_{L_\lambda^p(0,T)}.
\end{split}
\end{equation*}
Similarly, it follows from the condition $\|\partial_t^\alpha m_1\|_{C[0,T]}\le c_1$ that
\begin{align*}
\|   \mathrm{I}_2 \|_{L_\lambda^p(0,T)} & \le \Big\| \frac{ m_1  \partial_t^\alpha (m_2  -  m_1) }{m_1m_2} \Big\|_{L_\lambda^p(0,T)}
+ \Big\| \frac{ (m_1 - m_2) \partial_t^\alpha m_1 }{m_1m_2} \Big\|_{L_\lambda^p(0,T)} \\
&\le c \| \partial_t^\alpha (m_1  -  m_2)  \|_{L_\lambda^p(0,T)} + c \|  m_1 - m_2 \|_{L_\lambda^p(0,T)},\\
\|   \mathrm{I}_3 \|_{L_\lambda^p(0,T)} & {\le \Big\| \frac{ (m_1 - m_2) }{m_1m_2} \Big\|_{L_\lambda^p(0,T)}\|g\|_{C([0,T];L^2(\partial\Omega))}
\le c \|  m_1 - m_2 \|_{L_\lambda^p(0,T)}.}
\end{align*}
In sum, the last three estimates together yield
\begin{equation}\label{eqn:q1-q2}
\begin{split}
\|  q_1 - q_2\|_{L_\lambda^p(0,T)} &  \le c \big(\| u(q_1) - u(q_2) \|_{L_\lambda^p(0,T;L^2\II)} + \| \partial_t^\alpha (m_1  -  m_2)  \|_{L_\lambda^p(0,T)} + \|  m_1 - m_2 \|_{L_\lambda^p(0,T)}\big),
\end{split}
\end{equation}
{where the constant $c$ depends on $m_*$, $m^*$, $ \|\partial_t^\alpha m_1\|_{C[0,T]}$,   $\|g\|_{C([0,T];L^2(\partial\Omega))}$, and Lipschitz constant of $f$, but not on $\lambda$ and $p$.} Next, we  bound $\| u(q_1) - u(q_2) \|_{L_\lambda^p(0,T;L^2\II)}$. Note that $w = u(q_1) - u(q_2)$ satisfies
$$ \partial_t^\alpha w(t) + A w(t) = u(q_1)(q_2(t) - q_1(t)) + (f(u(q_1),t) - f(u(q_2),t)) - q_2(t) w(t),\quad t>0 \,\,\text{with} ~~ w(0)=0.  $$
Using the representation \eqref{eqn:sol}, $w$ can be represented by
\begin{align*}
w(t) &= \int_0^t E(t-s) u(s;q_1)(q_2(s) - q_1(s))\, \ds + \int_0^t E(t-s) (f(u(s;q_1),s) - f(u(s;q_2),s)) \,\ds \\
&\quad -  \int_0^t E(t-s) w(s) q_2(s)\,\ds = \sum_{i=1}^3 \mathrm{II}_i.
\end{align*}
We bound the three terms $\mathrm{II}_i$ separately. From Lemma \ref{lem:op} (ii) and the estimate $\| u(q_1)  \|_{C([0,T];L^2\II)} \le c$ (cf. Assumption \ref{ass:stab} and Theorem \ref{thm:sol-reg}), we deduce
\begin{align*}
e^{-\lambda t} \| \mathrm{II}_1(t)\|_{L^2\II} &\le  e^{-\lambda t}\int_0^t  \| E(t-s)\|_{L^2\II\rightarrow L^2\II} \| u(s;q_1) \|_{L^2\II} |q_2(s) - q_1(s)| \, \ds \\
&\le ce^{-\lambda t}\int_0^t  (t-s)^{\alpha-1}  |q_2(s) - q_1(s)| \, \ds.
\end{align*}
Then Young's inequality for convolution leads to
\begin{align*}
 \| \mathrm{II}_1 \|_{L_\lambda^p(0,T;L^2\II)} &\le   c \Big(\int_0^T\Big|\int_0^t   e^{-\lambda (t-s)}(t-s)^{\alpha-1}  |q_1(s) - q_2(s)|e^{-\lambda s} \, \ds\Big|^p \,\dt\Big)^{\frac1p}\\
 &\le c \int_0^T e^{-\lambda t}  t^{\alpha -1} \,\dt \,  \Big(\int_0^T  (|q_1(t) - q_2(t)|e^{-\lambda t})^p \, \dt\Big)^\frac1p \\
 & =  c \int_0^T e^{-\lambda t}  t^{\alpha -1} \,\dt \,  \| q_1 -q_2  \|_{L_\lambda^p(0,T)}.
\end{align*}
Then direct computation with the definition of the Gamma function $\Gamma(z)$ yields
\begin{align}\label{eqn:weight-est}
  \int_0^T e^{-\lambda t}  t^{\alpha -1} \,\dt = \lambda^{-\alpha}\int_0^{\lambda T} e^{-y}  y^{\alpha -1} \,\dy \le  \lambda^{-\alpha}\int_0^{\infty} e^{-y}  y^{\alpha -1} \,\dy = \Gamma(\alpha) \lambda^{-\alpha},
\end{align}
and hence we arrive at
\begin{align*}
 \| \mathrm{II}_1 \|_{L_\lambda^p(0,T;L^2\II)} &\le  c \lambda^{-\alpha}  \| q_1 -q_2  \|_{L_\lambda^p(0,T)}.
\end{align*}
For the term $\mathrm{II}_2$, from Lipschitz stability of $f$ in $u$ and Lemma \ref{lem:op} (ii), we derive
\begin{align*}
  \| \mathrm{II}_2(t)\|_{L^2\II} &\le   \int_0^t  \| E(t-s)\|_{L^2\II\rightarrow L^2\II} \| w(s) \|_{L^2\II} \, \ds \le  c\int_0^t (t-s)^{\alpha-1} \| w(s) \|_{L^2\II} \, \ds .
\end{align*}
Then repeating the argument for the term $\mathrm{II}_1$ leads to
\begin{align*}
  \| \mathrm{II}_2 \|_{L_\lambda^p(0,T;L^2\II)} &\le  c \lambda^{-\alpha}  \| w \|_{L_\lambda^p(0,T;L^2\II)}.
\end{align*}
Since $q_2 \in \mathcal{Q}$, $ \mathrm{II}_3$ can be bounded similarly as $\| \mathrm{II}_3 \|_{L_\lambda^p(0,T;L^2\II)} \le  c\lambda^{-\alpha}  \| w \|_{L_\lambda^p(0,T;L^2\II)}$. Hence, we arrive at
\begin{align*}
\| w \|_{L_\lambda^p(0,T;L^2\II)} \le  c \lambda^{-\alpha}\| w \|_{L_\lambda^p(0,T;L^2\II)} + c \lambda^{-\alpha}  \| q_1 -q_2  \|_{L_\lambda^p(0,T)},
\end{align*}
{with $c$ dependent of $\| u(q_1) \|_{L^\infty(0,T;L^2\II)}$ and Lipschitz constant of $f$, but not $\lambda$.}
Then choosing a sufficiently large $\lambda$ yields
\begin{align}\label{eqn:u1-u2}
\| u(q_1) - u(q_2) \|_{L_\lambda^p(0,T;L^2\II)} \le   c \lambda^{-\alpha}  \| q_1 -q_2  \|_{L_\lambda^p(0,T)}.
\end{align}
This together with \eqref{eqn:q1-q2} leads to
\begin{equation}\label{eqn:q1-q2-2}
\begin{split}
\|  q_1 - q_2\|_{L_\lambda^p(0,T)} &  \le c \big(\lambda^{-\alpha}  \| q_1 -q_2  \|_{L_\lambda^p(0,T)} + \| \partial_t^\alpha (m_1  -  m_2)  \|_{L_\lambda^p(0,T)} + \|  m_1 - m_2 \|_{L_\lambda^p(0,T)}\big).
\end{split}
\end{equation}
Again, we choose a sufficiently large $\lambda$ and obtain
\begin{equation*}
\|  q_1 - q_2\|_{L_\lambda^p(0,T)}  \le c \big(  \| \partial_t^\alpha (m_1  -  m_2)  \|_{L_\lambda^p(0,T)} + \|  m_1 - m_2 \|_{L_\lambda^p(0,T)}\big).
\end{equation*}
By the equivalence between the norms $\|\cdot \|_{L_\lambda^p(0,T)}$ and $\|\cdot \|_{L^p(0,T)}$, we obtain
\begin{equation}\label{eqn:stab-00}
\| q_1 - q_2  \|_{L^p(0,T)} \le c\big( \| m_1  - m_2   \|_{L^p(0,T)} + \| \partial_t^\alpha(m_1  - m_2)   \|_{L^p(0,T)} \big).
\end{equation}
Moreover, since $m_1(0) = m_2(0)$, by  \cite[Theorem 2.13(ii), p. 45]{Jin:2021}, we have
$(m_1 - m_2) (t)= \tfrac{1}{\Gamma(\alpha)}\int_0^t(t-s)^{\alpha-1} \partial_s^\alpha(m_1 - m_2)(s) \,\ds$.
This identity and Young's inequality for convolution lead to
$$ \| m_1 - m_2 \|_{L^p(0,T)}
\le c \int_0^T t^{\alpha-1}\,\dt \| \partial_t^\alpha(m_1 - m_2) \|_{L^p(0,T)}\le c\| \partial_t^\alpha(m_1 - m_2) \|_{L^p(0,T)}.$$
This and \eqref{eqn:stab-00} imply the desired stability result in the theorem.
\end{proof}

\begin{remark}
Theorem \ref{thm:stab} shows conditional Lipschitz stability for the inverse potential problem. This is expected, since the unknown potential $q(t)$ and the observational data $m(t)$ are aligned in direction, following the folklore theorem \cite{JinRundell:2015}. The Lipschitz stability estimate indicates that IPP amounts to an $\alpha$th order derivative loss, and hence as the fractional order $\alpha$ increases, IPP becomes more ill-posed. The estimate guarantees stable numerical reconstruction, which is also confirmed by the error estimate for the reconstruction algorithm below. This result is largely comparable with that for the observation at one interior point \cite[Theorem 2.2]{FujishiroKian:2016}.
\end{remark}

The stability estimate in Theorem \ref{thm:stab} naturally motivates developing a reconstruction algorithm to recover the time-dependent potential $q$ from the measurement $m$ with error estimates. Next we present a simple iterative algorithm and show the linear convergence in the weighted $L^p(0,T)$ norm.
\begin{proposition}\label{cor:iter}
Let Assumption \ref{ass:stab} hold. Then for any initial guess $q^0 \in \mathcal{Q}$, consider the following iteration
\begin{equation*}
    q^{k+1}(t) = P_{\mathcal{Q}}\Big[\frac{\int_{\Omega}f(u(t;q^{k}),t) \dx - \partial^{\alpha}_t m (t) +\int_{\partial \Omega} g(t) \ds  }{m(t)}\Big],
\end{equation*}
where $P_{\mathcal{Q}}$ denotes a cut-off operation such that for all $q\in C[0,T]$,
\begin{equation}
P_{\mathcal{Q}} q(t) = \min(\max(q(t),0),c_0),\quad \forall t\in[0,T].
\end{equation}
The iteration converges to $q^\dag$ in $L_\lambda^p(0,T)$ for sufficiently large $\lambda$ in the sense that
\begin{equation}\label{eqn:conv}
  \| q^{k}  -  q^\dag \|_{L_\lambda^p(0,T)} \le  (c\lambda^{-\alpha})^{k} \| q^0 - q^\dag  \|_{L_\lambda^p(0,T)}, \quad k=1,2,\ldots,
\end{equation}
{where the constant $c$ depends on $c_0$, $\|  u(q^\dag)\|_{C([0,T];L^2\II)}$ and Lipschitz constant of $f$, but not on $k$.}
\end{proposition}

\begin{proof}
We define an operator $K: \mathcal{Q} \rightarrow \mathcal{Q}$ such that for any $q\in\mathcal{Q}$,
$$ K q = P_{\mathcal{Q}}\Big[\frac{\int_{\Omega}f(u(t;q),t) \dx - \partial^{\alpha}_t m (t) +\int_{\partial \Omega} g(t) \ds  }{m(t)}\Big]. $$
One can show that the operator $K$ is indeed well-defined by repeating the argument of Theorem \ref{thm:sol-reg} and the fact that $ P_{\mathcal{Q}} q \in \mathcal{Q}$ for any $q\in C[0,T]$.
Moreover, by definition, $q^\dag$ is a fixed point of $K$ in $\mathcal{Q}$.
For any $k=0,1,\ldots$, using the stability of $P_{\mathcal{Q}}$, i.e., $|P_{\mathcal{Q}} q(t) - q^\dag(t)|\le |q(t) - q^\dag(t)|$ for all $q\in C[0,T]$ and $t\in[0,T]$, we obtain
\begin{align*}
 |q^{k+1}(t) - q^\dag(t)|
 &=  \Big|P_{\mathcal{Q}}\Big[\frac{\int_{\Omega}f(u(t;q^{k}),t) \dx - \partial^{\alpha}_t m (t) +\int_{\partial \Omega} g(t) \ds  }{m(t)}\Big] - q^\dag(t)\Big|\\
 &\le \Big| \frac{\int_{\Omega}f(u(t;q^{k}),t) \dx - \partial^{\alpha}_t m (t) +\int_{\partial \Omega} g(t) \ds  }{m(t)}  - q^\dag(t)\Big| \\
 &= \Big| \frac{\int_{\Omega}f(u(t;q^{k}),t) - f(u(t;q^\dag),t) \dx   }{m(t)}  \Big|.
\end{align*}
By the Lipschitz continuity of $f$ in $u$ and Assumption \ref{ass:stab} (i), we have
\begin{align*}
 |q^{k+1}(t) - q^\dag(t)|
  \le c   \int_{\Omega} |f(u(t;q^{k}),t) - f(u(t;q^\dag),t)| \dx    \le c \| u(t;q^{k}) - u(t;q^\dag) \|_{L^2\II},
\end{align*}
{with $c$ dependent of $m_*$ and Lipschitz constant of $f$.}
Then taking the $L^p_\lambda(0,T)$ norm on both sides gives
\begin{align*}
 \|q^{k+1}  - q^\dag \|_{L_\lambda^p(0,T)}   \le c \| u(t;q^{k}) - u(t;q^\dag) \|_{L_\lambda^p(0,T;L^2\II)}.
\end{align*}
Now we use the estimate \eqref{eqn:u1-u2} with   $q_1 = q^\dag$ and $q_2 = q^k$, and deduce
\begin{align*}
 \|q^{k+1}  - q^\dag \|_{L_\lambda^p(0,T)}    \le c \lambda^{-\alpha} \|  q^{k}  -   q^\dag  \|_{L_\lambda^p(0,T)} ,
\end{align*}
{where the constant $c$ is independent of $k$ and $\lambda$.} Then for sufficiently large $\lambda$ such that $c \lambda^{-\alpha}< 1$, we conclude that
the sequence $\{q^k\}_{k=0}^\infty$ converges to $q^\dag$ linearly in the sense of \eqref{eqn:conv}.
\end{proof}

\section{Numerical scheme and error estimate}\label{sec:recon}
In practice, we have to discretize the direct problem \eqref{eqn:fde} which incurs additional discretization errors, in addition to data noise. In this section, we shall develop a fully discrete numerical scheme for solving \IPP.
To this end, we introduce a fully discrete scheme based on backward Euler convolution quadrature in time \cite[Chapter 3]{JinZhou:2023book} and Galerkin finite element method in space \cite[Chapter 2]{JinZhou:2023book}. Then we present a reconstruction algorithm to recover the potential $q$ from the noisy observational
data $ m_\delta $.  Finally, we establish an \textit{a priori} error bound which provides guidelines to choose the
(space / time) mesh sizes $h$ and $\tau$ according to the noise level $\delta$.

\subsection{Numerical scheme for solving the direct problem}\label{sec:fully}
The literature on the numerical approximation of the  subdiffusion
model is vast, see e.g., \cite{JinLazarovZhou:2019, JinZhou:2023book} for recent overviews of existing schemes.
In this work, we employ convolution quadrature (CQ) to discretize the fractional derivative $\partial_t^\alpha u$ on uniform grids (cf. \cite{Lubich:1986,Lubich:1988} and \cite[Chapter 3]{JinZhou:2023book}).
Let $\{t_n=n\tau\}_{n=0}^N$ be a uniform partition of the interval $[0,T]$,
with a time step size $\tau=T/N$. Then the time stepping scheme for problem \eqref{eqn:fde}
reads: given $u^0(q)=u_0$, find $u^n(q) \in H^1\II$, $n=1,2,\ldots,N$, such that
\begin{align}\label{eqn:step-0}
   (\bar \partial_\tau^\alpha u^n(q), v) + (a\nabla u^n(q), \nabla v) + q(t_n)(u^n(q),v) = (f(u^{n-1}(q),t_n),v) + (g(t),v)_{L^2(\partial\Omega)}, \quad\forall v\in H^1\II.
\end{align}
The notation $\bar\partial_\tau^\alpha \varphi^n$ denotes the
backward Euler CQ approximation of $\partial_t^\alpha \varphi(t_n) $ (with $\varphi^j=\varphi(t_j)$):
\begin{equation}\label{eqn:CQ-BE}
  \bar\partial_\tau^\alpha \varphi^n = \tau^{-\alpha} \sum_{j=0}^n \omega_j^{(\alpha)} (\varphi^{n-j} - \varphi^0) ,\quad\mbox{ with } (1-\xi)^\alpha=\sum_{j=0}^\infty \omega_j^{(\alpha)}\xi^j.
\end{equation}
Note that the weights $\omega_j^{(\alpha)}$ are given explicitly by
$\omega_j^{(\alpha)} = (-1)^j\frac{\Gamma(\alpha+1)}{\Gamma(\alpha-j+1)\Gamma(j+1)}$.

For the space discretization, we employ the standard Galerkin finite element method \cite[Chapter 2]{JinZhou:2023book}.
We divide the domain $\Omega$ into a quasi-uniform simplicial triangulation $\mathcal{T}_h$ with a mesh size $h$.
Over the triangulation $\mathcal{T}_h$, we define a conforming piecewise linear finite element space $X_h\subset H^1(\Omega)$ by
\begin{equation*}
X_{h}:=\{v_{h}\in H^1(\Omega):	v_{h}|_{T}\in P_1(K) ,\,\,\, \forall K\in\mathcal{T}_h\},
\end{equation*}
where $P_1(K)$ denotes the set of linear polynomials on the element $K$.
On the finite element space $X_{h}$,
we define the standard $L^2(\Omega)$-projection $P_h:L^2(\Omega)\to X_h$ by
\begin{equation*}
	(P_hv,\varphi_h)=(v,\varphi_h), \quad \forall v\in L^2(\Omega),\  \varphi_h\in X_h,
\end{equation*}
and Ritz projection $R_h: H^1\II \to X_h$ such that
\begin{equation*}
	(a \nabla R_hv,\nabla \varphi_h)=(a \nabla v,\nabla \varphi_h)~~\quad\text{and}~~\int_\Omega R_hv \, \dx = \int_\Omega  v \,\dx , \quad \forall v\in H^1(\Omega),\  \varphi_h\in X_h.
\end{equation*}
Then the following approximation result holds:
\begin{equation}\label{inequ: L2 proj approx}
	\|v-P_hv\|_{L^2(\Omega)} +  \|v-R_hv\|_{L^2(\Omega)}\leq  ch^{2}\|v\|_{H^2(\Omega)}, \quad \forall v\in H^{2}(\Omega).
\end{equation}
The fully discrete scheme for problem \eqref{eqn:fde}
reads: given $u_h^0(q)= P_h u_0$, find $u_h^n(q) \in X_h$, $n=1,\ldots,N$, such that
\begin{equation}\label{eqn:fully}
\begin{aligned}
   (\bar \partial_\tau^\alpha u_h^n(q), \varphi_h) + (a\nabla u_h^n(q),& \nabla \varphi_h) + q(t_n)(u_h^n(q),\varphi_h) \\
     &= (f(u_h^{n-1}(q),t_n),\varphi_h) + (g(t),\varphi_h)_{L^2(\partial\Omega)},\quad \forall \varphi_h\in X_h.
\end{aligned}
\end{equation}

To analyze the scheme \eqref{eqn:fully}, we need preliminary estimates for the linear problem. Consider the linear problem
(with $q^\dag \equiv 0$): find $v(t) \in H^1\II$ such that
\begin{equation}\label{eqn:linear}
(\partial_t^\alpha v(t), \varphi) + (a\nabla v(t), \nabla \varphi) = (f(t),\varphi) + (g(t),\varphi)_{L^2(\partial\Omega)}, \quad \forall \varphi\in H^1(\Omega), \forall t\in(0,T],\quad \text{with} ~~v(0) = v_0,
\end{equation}
and its fully discrete scheme: given $v_h^0 = R_h v_0$, find $v_h^n\in X_h$ for $n=1,2,\ldots,N$ such that
\begin{align}\label{eqn:linear-fully}
   (\bar \partial_\tau^\alpha v_h^n , \varphi_h) + (a\nabla v_h^n , \nabla \varphi_h)  = (f(t_n),\varphi_h) + (g(t_n),\varphi_h)_{L^2(\partial\Omega)}, \quad \forall \varphi_h\in X_h.
\end{align}
\begin{lemma}\label{lem:error-linear}
Let $v$ and $v_h^n$ solve problems \eqref{eqn:linear} and \eqref{eqn:linear-fully}, respectively, with $v_0\in H^2(\Omega)$, and $f \in W^{1,1}(0,T;L^2(\Omega))$. Then under Assumption \ref{ass:stab}, the following error estimate holds
\begin{equation*}
   \|v(t_n)- v_h^n\|_{L^2(\Omega)} \leq c(h^2+t_n^{\alpha-1}\tau),
\end{equation*}
{where the constant $c$ depends on $\|g\|_{C^2([0,T];H^\frac12(\partial\Omega))}$, $\|v_0\|_{H^2(\Omega)}$ and $\|f\|_{W^{1,1}(0,T;L^2(\Omega))}$, but not on $h$, $\tau$ and $n$.}
\end{lemma}
\begin{proof}
The analysis is complicated by the presence of a nonzero Neumann data. First, we define a Neumann operator $\mathcal{N}:H^{-\frac12}(\partial\Omega) \rightarrow H^{1}\II$ such that for any $\psi \in H^{-\frac12}(\partial\Omega)$
$$ (a \nabla \mathcal{N} \psi, \nabla \varphi) = (\psi,\varphi)_{\partial\Omega}~~\text{and}~~\int_\Omega \mathcal{N} \psi \,\dx = 0,\quad \forall \varphi\in H^{1}\II.  $$
Then the function $w = v- \mathcal{N}g$ satisfies
\begin{equation}\label{eqn:linear-w}
\left\{\begin{aligned}
    (\partial_t^\alpha w(t), \varphi) + (a\nabla w(t), \nabla \varphi) &= (f(t) - \mathcal{N} \partial_t^\alpha g(t),\varphi), \quad \forall \varphi\in H^1(\Omega), \forall t\in(0,T],\\
    w(0) &= v_0 - \mathcal{N}g(0).
\end{aligned}\right.
\end{equation}
By construction, $w$ satisfies a zero Neumann boundary condition. Hence, by  \eqref{eqn:sol}, we have
$$ w(t) = F(t)(v_0 - \mathcal{N}g(0)) + \int_0^t E(t-s) [f(s) - \mathcal{N} \partial_s^\alpha g(s)]\,\ds. $$
Consequently,
$$v(t) = \mathcal{N}g(t) + F(t)(v_0 - \mathcal{N}g(0)) + \int_0^t E(t-s) [f(s) - \mathcal{N} \partial_s^\alpha g(s)]\,\ds. $$
Now consider the spatially semidiscrete scheme: find $v_h(t) \in X_h$ with $v_h(0) = R_h v_0$ such that
\begin{equation}\label{eqn:linear-fem}
(\partial_t^\alpha v_h(t), \varphi_h) + (a\nabla v_h(t), \nabla \varphi_h) = (f(t),\varphi_h) + (g(t),\varphi_h)_{L^2(\partial\Omega)},\quad \forall \varphi_h\in X_h,  \forall t\in(0,T].
\end{equation}
We define the discrete Neumann operator $\mathcal{N}_h: H^{-\frac12}(\partial\Omega)\rightarrow X_h$ by
$$ (a \nabla \mathcal{N}_h \psi, \nabla \varphi_h) = (\psi,\varphi_h)_{\partial\Omega}~~\text{and}~~\int_\Omega \mathcal{N}_h \psi \,\dx = 0,\quad \forall \varphi\in X_h.$$
Then the function $w_h = v_h - \mathcal{N}_hg$ satisfies
\begin{equation*}
\left\{\begin{aligned}
    (\partial_t^\alpha w_h(t), \varphi_h) + (a\nabla w_h(t), \nabla \varphi_h) &= (f(t) - \mathcal{N}_h \partial_t^\alpha g(t),\varphi_h),\quad \forall \varphi_h\in X_h,  ~~\forall t\in(0,T],\\ w_h(0) &= R_hv_0 - \mathcal{N}_h g(0).
\end{aligned}\right.
\end{equation*}
Similarly, the solution $v_h$ can be represented by
\begin{equation*}
\begin{split}
 v_h(t) = \mathcal{N}_h g(t) + F_h(t) (R_h v_0 - \mathcal{N}_h g(0)) + \int_0^t E_h(t-s) [P_h f(s) - \mathcal{N}_h \partial_s^\alpha g(s)]\,\ds,
\end{split}
\end{equation*}
where $F_h(t)$ and $E_h(t)$ are the spatially discrete analogues of $F(t)$ and $E(t)$ in \eqref{eqn:sol} \cite[Section 3.3]{JinZhou:2023book}
\begin{equation*}
\begin{split}
F_h(t):=\frac{1}{2\pi {\rm i}}\int_{\Gamma_{\theta,\delta }}e^{zt} z^{\alpha-1} (z^\alpha+A_h)^{-1}\, \dz \quad\mbox{and}\quad
E_h(t):=\frac{1}{2\pi {\rm i}}\int_{\Gamma_{\theta,\delta}}e^{zt}  (z^\alpha +A_h)^{-1}\, \dz ,
\end{split}
\end{equation*}
respectively, where the operator $A_h:X_h\rightarrow X_h$ is defined such that $(A_h\psi_h, \varphi_h) = (a\nabla \psi_h, \nabla \varphi_h)$ for all $ \psi_h, \varphi_h \in X_h$, with $\int_\Omega \psi_h \dx = 0$.
Since $R_h \mathcal{N} \psi = \mathcal{N}_h \psi$, the error $v_h-v$ is given by
\begin{equation*}
\begin{split}
v_h(t) - v(t) &= (R_h-I) \mathcal{N}g(t)  + (F_h(t) R_h - F(t))( v_0 - \mathcal{N} g(0))   \\
 &\quad + \int_0^t E_h(t-s) (P_h - R_h) \mathcal{N}\partial_s^\alpha g(s) \,\ds\\
 &\quad +  \int_0^t E_h(t-s)   P_h (f(s) - \mathcal{N}\partial_s^\alpha g(s)) - E(t-s)(f(s) - \mathcal{N}\partial_s^\alpha g(s)) \,\ds =: \sum_{j=1}^4 \mathrm{I}_j.
 \end{split}
\end{equation*}
The  approximation property in \eqref{inequ: L2 proj approx}
immediately implies
\begin{equation*}
\begin{split}
 \|  \mathrm{I}_1 + \mathrm{I}_3 \|_{L^2\II} &\le c h^2\Big( \| \mathcal{N}g(t) \|_{H^2(\Omega)} + \int_0^t (t-s)^{\alpha-1}  \| \mathcal{N}\partial_s^\alpha g(s)\|_{H^2\II}\,\ds \Big)\\
 & \le c h^2\Big( \| g(t) \|_{H^{\frac12}(\partial\Omega)} + \int_0^t (t-s)^{\alpha-1}  \| \partial_s^\alpha g(s)\|_{H^{\frac12}(\partial\Omega)} \,\ds\Big) \\
 &\le c h^2 \| g \|_{C^1([0,T];H^{\frac12}(\partial\Omega))}.
 \end{split}
\end{equation*}
Meanwhile, by applying the argument in \cite[Theorems 2.4 and 2.5 (ii)]{JinZhou:2023book}, we obtain
\begin{equation*}
\begin{split}
 \|  \mathrm{I}_2 + \mathrm{I}_4 \|_{L^2\II} &\le c h^2 \Big(\| v_0 - \mathcal{N} g(0) \|_{H^2\II} + \| f(0) - \mathcal{N} (\partial_t^{\alpha}g)(0) \|_{L^2\II} + \int_0^t \|\partial_s(f(s) - \mathcal{N} \partial_t^{\alpha}g(s))\|_{L^2\II} \,\ds\Big) \\
 &\le ch^2 \Big(\| v_0 \|_{H^2\II} + \| f(0) \|_{L^2\II} + \| g \|_{C^2([0,T];H^\frac12(\partial\Omega))} + \int_0^t \| f'(s) \|_{L^2\II}\,\ds\Big) .
 \end{split}
\end{equation*}
In sum, we arrive at the following bound
\begin{equation}\label{eqn:err-fem-linear}
\|v_h(t) - v(t)   \|_{L^2\II} \le  ch^2 \Big(\| v_0 \|_{H^2\II} + \| f(0) \|_{L^2\II} + \| g \|_{C^2([0,T];H^\frac12(\partial\Omega))} + \int_0^t \| f'(s) \|_{L^2\II}\,\ds\Big)  .
\end{equation}
Upon letting $w_h^n = v_h^n - \mathcal{N}_h g(t_n)$, we have $w_h^0 = R_hv_0 - \mathcal{N}_h g(0) = R_h (v_0 - \mathcal{N}g(0))$ and for $n=1,2,\ldots,N$,
\begin{align}\label{eqn:linear-fully-w}
   (\bar \partial_\tau^\alpha w_h^n , \varphi_h) + (a\nabla w_h^n , \nabla \varphi_h)  = (f(t_n) - \mathcal{N}_h \bar \partial_\tau^\alpha g(t_n),\varphi_h),\quad \forall  \varphi_h\in X_h.
\end{align}
The solution $w_h^n \in X_h$ can be represented by
$$ w_h^n = F_{h,\tau}^n w_h^0 + \tau \sum_{j=1}^n E_{h,\tau}^{n-j}  (P_hf(t_j) - \mathcal{N}_h \bar \partial_\tau^\alpha g(t_j)),$$
where the fully discrete solution operators $F_{h, \tau}^n $ and $E_{h, \tau}^n $ are defined respectively by (see e.g., \cite[equations (4.3)-(4.5)]{JinLiZhou:2019} and \cite[equation (9.24)]{JinZhou:2023book})
\begin{equation}\label{eqn:FEht-0}
\begin{aligned}
F_{h, \tau}^n  &= \frac{1}{2\pi\mathrm{i}}\int_{\Gamma_{\theta,\sigma}^\tau } e^{zt_n} {e^{-z\tau}} \delta_\tau(e^{-z\tau})^{\alpha-1}({ \delta_\tau(e^{-z\tau})^\alpha}+{A_h)^{-1}}\,\dz,\\
E_{h, \tau}^n  &= \frac{1}{2\pi\mathrm{i}}\int_{\Gamma_{\theta,\sigma}^\tau } e^{zt_n} ({ \delta_\tau(e^{-z\tau})^\alpha}+{A_h)^{-1}}\,\dz ,
\end{aligned}
\end{equation}
with $\delta_\tau(\xi)=(1-\xi)/\tau$ and the contour
$\Gamma_{\theta,\sigma}^\tau :=\{ z\in \Gamma_{\theta,\sigma}:|\Im(z)|\le {\pi}/{\tau} \}$, with $\theta\in(\pi/2,\pi)$ being close to $\pi/2$
(oriented with an increasing imaginary part).
Next we define an auxiliary function $\bar w_h^n$ such that $\bar w_h^0 = w_h^0$ and for $n=1,2,\ldots,N$,
\begin{align*}
   (\bar \partial_\tau^\alpha \bar w_h^n ,\varphi_h) + (a\nabla \bar w_h^n , \nabla \varphi_h)  = (f(t_n) - \mathcal{N}_h  \partial_t^\alpha g(t_n),\varphi_h), \quad \forall \varphi_h \in X_h.
\end{align*}
Note that $\bar \theta_h^n = w_h^n - \bar w_h^n$ satisfies $\bar \theta_h^n = 0$ and for $n=1,2,\ldots,N$,
\begin{align*}
   (\bar \partial_\tau^\alpha \bar \theta_h^n ,\varphi_h) + (a\nabla  \bar \theta_h^n , \nabla \varphi_h)
   = (   (\partial_t^\alpha -\bar \partial_\tau^\alpha) \mathcal{N}_h g(t_n),\varphi_h),\quad \forall \varphi_h \in X_h.
\end{align*}
Recall the following stability estimate  \cite[Lemma 3.4]{WuZhou:2021}
\begin{equation}\label{eqn:Ehn-stab}
\| E_{h,\tau}^{n-j}  \|_{L^2\II\rightarrow L^2\II} \le c(t_{n+1} - t_j)^{\alpha-1}.
\end{equation}
Then we can derive
\begin{align*}
   \| \bar \theta_h^n \|_{L^2\II}
&\le c \tau \sum_{j=1}^n \| E_{h,\tau}^{n-j}  \|_{L^2\II\rightarrow L^2\II} \| (\partial_t^\alpha -\bar \partial_\tau^\alpha) \mathcal{N}_h g(t_n)\|_{L^2\II} \\
&\le c \tau \sum_{j=1}^n (t_{n+1} - t_j)^{\alpha-1} \| (\partial_t^\alpha -\bar \partial_\tau^\alpha) \mathcal{N}_h g(t_n)\|_{L^2\II} \le c\tau \|   g \|_{C^2(0,T;L^2\II)},
\end{align*}
where the last inequality follows from the   truncation error estimate \cite[Theorem 3.1 with $k=1$]{JinZhou:2023book}
\begin{align*}
 \| (\partial_t^\alpha - \bar \partial_\tau^\alpha)\mathcal{N}_h g(t_n) \|_{L^2\II}
 &\le c\tau  \| \mathcal{N}_h g \|_{C^2(0,T;L^2\II)}
= c\tau  \| R_h \mathcal{N} g \|_{C^2(0,T;L^2\II)}\\
 & \le c\tau   \|  \mathcal{N} g \|_{C^2(0,T;H^1\II)} \le c\tau  \|   g \|_{C^2(0,T;L^2({\partial\Omega)})}.
\end{align*}
Moreover, using the estimate \cite[Theorem 3.4]{JinZhou:2023book}, we obtain
\begin{align*}
 \| w_h(t_n) - w_h^n \|_{L^2\II}
 &\le c \tau \big(t_n^{\alpha-1}\| P_hf(0) - R_h \mathcal{N}   \partial_t^\alpha g(0) - A_h w_h^0 \|_{L^2\II}\big) \\
 &\quad +c\tau \int_0^t (t_n-s)^{\alpha-1} \| P_hf'(s) - \mathcal{N}_h  \partial_s \partial_s^\alpha g(s) \|_{L^2\II}\\
 &\le c \tau \big(t_n^{\alpha-1} (\| f(0)\|_{L^2\II} +  \| v_0 \|_{H^2\II}) + \| g \|_{C^2([0,T];L^2(\partial\Omega))}\big) \\
 &\quad + c\tau\int_0^t (t_n-s)^{\alpha-1} \| f'(s)\|_{L^2\II}   \,\ds.
\end{align*}
Therefore, the following estimate holds
\begin{align*}
 &\| v_h(t_n) - v_h^n \|_{L^2\II} = \| w_h(t_n) - w_h^n \|_{L^2\II}\\
  \le& c \tau \Big(t_n^{\alpha-1} (\| f(0)\|_{L^2\II} +  \| v_0 \|_{H^2\II}) + \| g \|_{C^2([0,T];L^2(\partial\Omega))} + \int_0^t (t_n-s)^{\alpha-1} \| f'(s)\|_{L^2\II}   \,\ds\Big).
\end{align*}
This together with the estimate \eqref{eqn:err-fem-linear} shows the desired error bound. This completes the proof of the lemma.
\end{proof}

The next lemma provides an error estimate for the fully discrete scheme \eqref{eqn:fully}.
\begin{lemma}\label{lem:fully-direct} Let $q^\dag\in \mathcal{Q}\cap C^1[0,T]$, $u_0 \in H^2\II$, $g\in C^2([0,T];H^{\frac12}(\partial\Omega))$ with $\partial_\nu u_0 = g(0)$ on $\partial\Omega$, and let $f(u,x,t)$ be smooth and globally Lipschitz.
Let $u = u(q^\dag)$ be the solution to  problem \eqref{eqn:fde} and $u_h^n = u_h^n(q^\dag)$ be the solution to
the numerical scheme  \eqref{eqn:fully}. Then there holds
$$\| u(t_n ) - u_h^n  \|_{L^2\II} \le c (\tau t_n^{\alpha-1} + h^2), $$
{where the constant $c$ depends on $\|q\|_{C^1[0,T]}$, $\|u_0\|_{H^2(\Omega)}$, $\|g\|_{C^2([0,T];H^\frac12(\partial\Omega))}$, Lipschitz constant of $f$ and terminal time $T$, but not on $h$, $\tau$ and $n$.}
\end{lemma}
\begin{proof}
We introduce an auxiliary function $\bar u_h^n$  with $\bar u_h^0 = R_h u_0$ and for  $n=1,2,\ldots,N$,
\begin{align}\label{eqn:fully-bar}
   (\bar \partial_\tau^\alpha \bar u_h^n , \varphi_h) + (a\nabla \bar u_h^n , \nabla \varphi_h)  = (f(u(t_n),t_n),\varphi_h)-q^\dag(t_n)(u(t_n),\varphi_h) + (g(t_n),\varphi_h)_{\partial\Omega}, \quad \forall \varphi_h\in X_h.
\end{align}
Then we split the error $e^n=u(t_n)-u_h^n$ into
$$ e^n = (u(t_n) - \bar u_h^n) + ( \bar u_h^n -  u_h^n) =: \vartheta^n + \varrho_h^n.$$
By Lemma \ref{lem:error-linear},
and the solution regularity in Theorem \ref{thm:sol-reg}, we deduce
\begin{equation}\label{eqn:err-nonlinear-1}
\begin{split}
 \| \vartheta^n \|_{L^2\II} \le c (\tau t_n^{\alpha-1} + h^2).
\end{split}
\end{equation}
Moreover, note that $\varrho_h^0 = 0$ and for $n=1,2,\ldots,N$,
\begin{align}\label{eqn:fully-uh-buh}
   (\bar \partial_\tau^\alpha \varrho_h^n , \varphi_h) + (a\nabla \varrho_h^n , \nabla \varphi_h)  = (f(u(t_n),t_n)-f(u_h^{n-1},t_n),\varphi_h)-q^\dag(t_n)(u(t_n)-u_h^{n},\varphi_h),\quad \forall \varphi_h \in X_h.
\end{align}
Then the error $\varrho_h^n $ can be represented using the operator $E_{h,\tau}^j$ in \eqref{eqn:FEht-0} as
$$ \varrho_h^n =  \tau \sum_{j=1}^n E_{h,\tau}^{n-j}  \big(P_h(f(u(t_j),t_j)-{f(u_h^{j-1},t_j)) -q^\dag(t_j)}(P_h u(t_j)-u_h^{j}) \big). $$
Using the \textsl{a priori} estimate \eqref{eqn:Ehn-stab} and the $L^2\II$ stability of $P_h$, we derive
\begin{equation*}
\begin{split}
 \| \varrho_h^n  \|_{L^2\II} &\le c \tau \sum_{j=1}^n (t_{n+1} - t_j )^{\alpha-1} \big(\|P_h(f(u(t_j),t_j)-f(u_h^{j-1},t_j))  \|_{L^2\II} + \| P_h u(t_j)-u_h^{j} \|_{L^2\II} \| q^\dag \|_{C[0,T]}\big)\\
&\le c \tau \sum_{j=1}^n (t_{n+1} - t_j )^{\alpha-1} \big(\| f(u(t_j),t_j)-f(u_h^{j-1},t_j)  \|_{L^2\II} + \| e^j \|_{L^2\II} \big).
\end{split}
\end{equation*}
Then by the triangle inequality and  Lipschitz continuity of $f$, we have
\begin{equation*}
\begin{split}
 \| \varrho_h^n  \|_{L^2\II} &\le  c \tau \sum_{j=1}^n (t_{n+1} - t_j )^{\alpha-1} \big(\| u(t_j) - u_h^{j-1}   \|_{L^2\II} + \| e^j \|_{L^2\II} \big) \\
&\le c \tau \sum_{j=1}^n (t_{n+1} - t_j )^{\alpha-1} \big(\| u(t_j) - u(t_{j-1})   \|_{L^2\II} + \| e^j \|_{L^2\II} +  \| e^{j-1} \|_{L^2\II}\big) \\
&\le c \tau t_n^{\alpha-1} + c \tau \sum_{j=1}^n (t_{n+1} - t_j )^{\alpha-1}\big(\| e^j \|_{L^2\II} +  \| e^{j-1} \|_{L^2\II}\big),
\end{split}
\end{equation*}{
where the last step follows from the estimate $\| u(t) \|_{L^2\II} + t^{1-\alpha} \| \partial_t u(t) \|_{L^2\II} \le c$ (cf. Theorem \ref{thm:sol-reg}) as
\begin{align*}
    &\tau \sum_{j=1}^n (t_{n+1} - t_j )^{\alpha-1} \| u(t_j) - u(t_{j-1})   \|_{L^2\II} \\
    \leq& \tau \sum_{j=1}^n (t_{n+1} - t_j )^{\alpha-1} \int_{t_{j-1}}^{t_j}\|\partial_s u(s)\|_{L^2\II} \ds\\
    \leq &c\tau \sum_{j=1}^n (t_{n+1} - t_j )^{\alpha-1} \int_{t_{j-1}}^{t_j}s^{\alpha-1} \ds\leq c\sum_{j=1}^n \int_{t_{j-1}}^{t_j}(t_{n} - s)^{\alpha-1}s^{\alpha-1}\ds \\
    =& c\int_0^{t_n}(t_n-s)^{\alpha-1}s^{\alpha-1}\ds
    \leq ct_n^{2\alpha-1} \le c_Tt_n^{\alpha-1}.
\end{align*}}
This together with \eqref{eqn:err-nonlinear-1} leads to
\begin{equation*}
 \| e^n  \|_{L^2\II} \le   c (h^2 + \tau t_n^{\alpha-1}) + c \tau \sum_{j=1}^n (t_{n+1} - t_j )^{\alpha-1}\big(\| e^j \|_{L^2\II} +  \| e^{j-1} \|_{L^2\II}\big).
\end{equation*}
Then applying the discrete Gronwall's inequality for convolution quadrature \cite[Theorem 10.2, p. 262]{JinZhou:2023book} gives the desired estimate, completing the proof of the lemma.
\end{proof}

\subsection{Numerical scheme for the inverse potential problem}\label{inverse_potential}
Now we discuss a numerical scheme for recovering the time-dependent potential $q(t)$.
Throughout, we assume that the measurement data  $m_\delta$ is noisy in the sense that
\begin{equation}\label{eqn:noise}
 \| m_\delta - m \|_{C[0,T]} = \delta \le \frac{m_*}{2}.
\end{equation}
Under Assumption \ref{ass:stab} (i),  $m_\delta$ is bounded away from zero such that
$$ 0 < \frac{m_*}{2} \le m_\delta(t) \le m^*+\frac{m_*}{2}, \quad \forall t\in [0,T]. $$
We will derive an error estimate for the numerical reconstruction in the discrete $\ell^p$ norm.
For $1\le p\leq \infty$, we denote by $\ell^p(X)$ the space of sequences $v^n\in X$, $n=1,\dots$,
such that\index{$\ell^p(X)$} $\|(v^n)_{n=1}^\infty\|_{\ell^p(X)}<\infty$, with
$$
\|(v^n)_{n=1}^\infty\|_{\ell^p(X)}:=
\left\{
\begin{aligned}
&\bigg(\sum_{n=1}^\infty\tau\|v^n\|_{X}^p\bigg)^{\frac{1}{p}},  &&\mbox{if}\,\,\, 1\le p<\infty,\\
&\sup_{n\ge 1}\|v^n\|_{X}, &&\mbox{if}\,\,\, p=\infty .
\end{aligned}\right.
$$
For a finite sequence $(v^n)_{n=1}^m\subset X$, we denote
$\|(v^n)_{n=1}^m\|_{\ell^p(X)}:=\|(v^n)_{n=1}^\infty\|_{\ell^p(X)}$,
by setting $v^n=0$ for $n>m$.

The following lemma provides a crucial estimate for the (fractional-order) numerical differentiation of the data $m_\delta(t)$. The estimate indicates the need of a proper selection of the time step size $\tau$ in order to optimally balance the data propagation error $\tau^{-\alpha}\delta$ and the time discretization error $\tau^{1/p}|\log \tau|$. This trade-off arises from the convolution quadrature approximation $\bar\partial_\tau^\alpha m_\delta(t_n)$ of the noisy data $m_\delta\in C[0,T]$. Note that if
$m_\delta$ has only the $C[0,T]$ regularity, the fractional derivative $\partial_t^\alpha m_\delta$ is actually ill-defined, and this is also the main source of computational challenges for IPP.

\begin{lemma}\label{lem:num-diff}
Let Assumption \ref{ass:stab}  hold and the data $m_\delta$ satisfy \eqref{eqn:noise}. Then for any $p\in{(1,\infty)}$, the convolution quadrature approximation $\bar\partial_\tau^\alpha m_\delta(t_n)$ defined by \eqref{eqn:CQ-BE} satisfies
$$  \|(\partial_t^\alpha m(t_n) - \bar\partial_\tau^\alpha m_\delta(t_n))_{n=1}^N\|_{\ell^p(\mathbb{R})}  \le c\big(   \delta \tau^{-\alpha} + \tau^{1/p}|\log\tau|\big),$$
{where the constant $c$ depends on $m$ and $p$, and not on $\delta$, $\tau$ and $n$.}
\end{lemma}
\begin{proof}
By the definition of convolution quadrature, for all $n=1,\ldots, N$,
\begin{align*}
| \bar\partial_\tau^\alpha m_\delta(t_n) - \bar\partial_\tau^\alpha m (t_n)| \le \tau^{-\alpha} \Big|\sum_{j=0}^n \omega_{n-j}^{(\alpha)}\big[(m_\delta(t_n)-m_\delta(0)) - (m(t_n)-m(0))\big]\Big|.
\end{align*}
It can be verified directly that $\omega_0^{(\alpha)}=1$,  $\omega_j^{(\alpha)}<0$ for $j\geq 1$
and $\sum_{j=0}^\infty \omega_j^{(\alpha)} = 0$. Hence, we have
\begin{align*}
| \bar\partial_\tau^\alpha m_\delta(t_n) - \bar\partial_\tau^\alpha m (t_n)| \le \tau^{-\alpha} \|m_\delta - m \|_{C[0,T]} \sum_{j=0}^n |\omega_{n-j}^{(\alpha)}|  \le c \delta \tau^{-\alpha}.
\end{align*}
Meanwhile, under Assumption \ref{ass:stab} (ii), we claim
\begin{align}\label{eqn:err-mt}
| \partial_t^\alpha m (t_n) - \bar\partial_\tau^\alpha m (t_n)| \le c \tau |\log \tau |  t_n^{-1}.
\end{align}
Thus we deduce
\begin{align*}
 \|[\partial_t^\alpha m(t_n) - \bar\partial_\tau^\alpha m_\delta(t_n)]_{n=1}^N\|_{\ell^p(\mathbb{R})}
& \le \|[\bar\partial_\tau^\alpha m_\delta(t_n) - \bar\partial_\tau^\alpha m (t_n)]_{n=1}^N\|_{\ell^p(\mathbb{R})} \\&\quad + \|[ \partial_t^\alpha m (t_n) - \bar\partial_\tau^\alpha m (t_n)]_{n=1}^N\|_{\ell^p(\mathbb{R})} \le c (\delta \tau^{-\alpha} + \tau^{1/p} |\log \tau|).
\end{align*}
Finally, we show the claim \eqref{eqn:err-mt}. Let
$w(t) = m(t) - m(0)$. Then \eqref{eqn:err-mt} is equivalent to
\begin{align}\label{eqn:err-mt-2}
| \partial_t^\alpha w (t_n) - \bar\partial_\tau^\alpha w (t_n)| \le \tau |\log \tau| t_n^{-1}.
\end{align}
Note that under Assumption \ref{ass:stab} (ii), i.e., $|\partial_t^\ell (\partial_t^\alpha w)(t)| \le c t^{-\ell}$ for $\ell=0,1$ and $|\partial_t^{s} w(t)| \le c t^{\alpha-s}$ for $s=1,2$, direct calculation with $c_w=(\partial_t^\alpha w)(0)$ yields
\begin{align*}
w(t) = [\omega * \partial_t (\partial_t^\alpha w)](t) + c_wp(t),\quad \text{with}~~ \omega(t) = \frac{t^\alpha}{\Gamma(1+\alpha)}.
\end{align*}
Indeed, by  \cite[Theorem 2.13(ii), p. 45]{Jin:2021}, since $w(0)=0$, and by integration by parts, there holds
\begin{align*}
    w(t) = (\omega'\ast \partial_t^\alpha w)(t) = (\omega\ast \partial_t(\partial_t^\alpha w)) - \omega(t-s)(\partial_s^\alpha w)|_{s=0}^{s=t} = (\omega\ast \partial_t\partial_t^\alpha w)+ c_w\omega(t).
\end{align*}
Then repeating the argument for \cite[equation (3.9)]{JinZhou:2023book} yields
\begin{equation}\label{eqn:err-dp}
    \left|\partial_t^\alpha \omega\left(t_n\right)-\bar{\partial}_\tau^\alpha \omega\left(t_{n}\right)\right| \leq c \tau t_{n+1}^{-1} .
\end{equation}
Next, the property of convolution implies
\begin{equation*}
\partial_t^\alpha w
= \partial_t^\alpha\big([\omega * \partial_t (\partial_t^\alpha w)](t) + \omega(t)\partial_t^\alpha w(0)\big)
   = (\partial_t^\alpha\omega)* (\partial_t (\partial_t^\alpha w))(t) + c_w\partial_t^\alpha \omega(t).
\end{equation*}
Let
$$
E_\tau(t)=\tau^{-\alpha} \sum_{n=0}^{\infty} \omega_n^{(\alpha)} \delta_{t_n}(t),
$$
with $\delta_{t_n}$ being the Dirac-delta function at $t_n$ (from the left side). Then we have
\begin{equation*}
\begin{aligned}
\bar{\partial}_\tau^\alpha w\left(t_n\right)
&=\left(E_\tau *\left(\omega* \partial_t (\partial_t^\alpha w)\right)\right)\left(t_n\right)+ c_w\bar\partial_\tau^\alpha \omega(t_n) \\
&=\left(\left(E_\tau *\omega\right) * \partial_t (\partial_t^\alpha w)\right)\left(t_n\right) + c_w\bar\partial_\tau^\alpha \omega(t_n).
\end{aligned}
\end{equation*}
Moreover, note that $\partial_t^\alpha \omega(t) = 1$. Then
for $n\ge2$ and $t\in(t_{n-1},t_n)$, we apply the estimate \eqref{eqn:err-dp} and obtain
\begin{equation}\label{eqn:err-dpt}
\left|\left(\partial_t^\alpha\omega-E_\tau *\omega\right)(t)\right| = \left|\left(\partial_t^\alpha \omega-E_\tau *\omega\right)(t_{n-1})\right| \le c \tau t_n^{-1} \le c \tau (t+\tau)^{-1}.
\end{equation}
Now we consider the following splitting
\begin{align*}
\left|\partial_t^\alpha w\left(t_n\right)-\bar{\partial}_\tau^\alpha w\left(t_n\right)\right|
\leq &\Big| \int_0^{\tau} (\partial_t^\alpha\omega - E_\tau*\omega)(t_n-s) (\partial_s \partial_s^\alpha) w(s) \,\ds\Big| \\
 & + \Big| \int_\tau^{t_n} (\partial_t^\alpha \omega - E_\tau*\omega)(t_n-s) (\partial_s \partial_s^\alpha) w(s) \,\ds\Big| \\
  &+ |c_w|\big|(\partial_t^\alpha - \bar{\partial}_\tau^\alpha ) \omega(t_n)\big| =: \sum_{i=1}^3 \mathrm{I}_i.
\end{align*}
Note that for $t\in (t_{n-1},t_n)$, since $\partial_t^\alpha \omega$ is constant valued and $E_\tau * p$ is piecewise constant, we have
$(\partial_t^\alpha \omega - E_\tau * \omega) (t) = \partial_t^\alpha \omega (t_{n-1})- (E_\tau * \omega)(t_{n-1}).$
This, \eqref{eqn:err-dp} and  Assumption \ref{ass:stab} (ii) imply
\begin{equation*}
\begin{aligned}
 \mathrm{I}_1
 &=  |\partial_t^\alpha \omega (t_{n-1})- (E_\tau * \omega)(t_{n-1})| \, \Big| \int_0^{\tau} (\partial_s\partial_s^\alpha) w(s) \,\ds\Big|\\
 &=  |\partial_t^\alpha \omega (t_{n-1})- (E_\tau * \omega)(t_{n-1})| \, \Big| \partial_t^\alpha w(\tau) - \partial_t^\alpha w(0)  \Big|\\
 &\le c|\partial_t^\alpha \omega (t_{n-1})- (E_\tau * \omega)(t_{n-1})| \le c \tau t_n^{-1}.
\end{aligned}
\end{equation*}
For the term $\mathrm{I}_2$, by \eqref{eqn:err-dpt}
and  Assumption \ref{ass:stab} (ii), we obtain
\begin{equation*}
\begin{aligned}
\mathrm{I}_2
\le c \tau \int_\tau^{t_n} (t_{n+1} - s)^{-1} |(\partial_s \partial_s^\alpha) w(s)| \,\ds
\le c \tau \int_\tau^{t_n} (t_{n+1} - s)^{-1} s^{-1}\,\ds
\le c \tau t_{n}^{-1} |\log \tau|.
\end{aligned}
\end{equation*}
Finally, the term $\mathrm{I}_3$ can be bounded using
\eqref{eqn:err-dp}
and  Assumption \ref{ass:stab} (ii) as
$$\mathrm{I}_3 \le c \tau t_{n}^{-1}
 \le c \tau t_{n}^{-1}|\log \tau|.$$
Combining the preceding bounds yields the claim \eqref{eqn:err-mt}.
\end{proof}

\begin{remark}
{Lemma \ref{lem:num-diff} discusses only the case $p\in(1,\infty)$, and does not cover the cases $p=1,\infty$. Indeed, the proof only gives
$\|(\partial_t^\alpha m(t_n) - \bar\partial_\tau^\alpha m_\delta(t_n))_{n=1}^N\|_{\ell^\infty(\mathbb{R})}  \le c\big(   \delta \tau^{-\alpha} + |\log\tau|\big),$
which does not tend to zero as $\tau\to0^+$. Also we have $\|(\partial_t^\alpha m(t_n) - \bar\partial_\tau^\alpha m_\delta(t_n))_{n=1}^N\|_{\ell^1(\mathbb{R})}  \le c\big(   \delta \tau^{-\alpha} + \tau |\log\tau|^2\big)$, involving an additional log factor.
}
\end{remark}

Next, we study a fully discrete iterative algorithm for recovering the time-dependent potential $q(t)$ at all discrete time levels.
To show the convergence, we employ the following weighted $\ell^p$ norm.
For $1\le p\leq \infty$ and a fixed $\lambda\ge0$, let  $\ell_\lambda^p(X)$ be the weighted norm for a sequences $v^n\in X$, $n=1,\dots$,
such that
$$
\|[v^n]_{n=1}^\infty\|_{\ell^p_\lambda(X)}:=
\left\{
\begin{aligned}
&\bigg(\tau\sum_{n=1}^\infty (e^{-\lambda t_n}\|v^n\|_{X})^p\bigg)^{\frac{1}{p}},  &&\mbox{if}\,\,\, 1\le p<\infty,\\
&\sup_{n\ge 1} e^{-\lambda t_n}\|v^n\|_{X}, &&\mbox{if}\,\,\, p=\infty .
\end{aligned}\right.
$$
Note that for any sequence $(v^n)_{n=1}^N$, there holds
\begin{equation}\label{eqn:quiv-disc}
\|[v^n]_{n=1}^N\|_{\ell^p_\lambda(X)} \le \|[v^n]_{n=1}^N\|_{\ell^p(X)} \le  e^{\lambda T}\|[v^n]_{n=1}^N\|_{\ell_\lambda^p(X)}.
\end{equation}
\begin{theorem}\label{thm:recon}
Let Assumption \ref{ass:stab} hold, and let the data $m_\delta$ satisfy \eqref{eqn:noise}.
Then for any initial guess $q_0=(q_0^n)_{n=1}^N \in {\mathcal{Q}_N}:=\{q=(q^n)_{n=1}^N: 0\le q^n \le c_0\}$,
consider the following fixed point iteration, which updates $q_{k+1}=(q_{k+1}^n)_{n=1}^N$ from  $q_{k}=(q_{k}^n)_{n=1}^N$ by
\begin{equation*}
    q_{k+1}^n = P_{\mathcal{A}_N}\Big[\frac{\int_{\Omega}f(u_h^n(q_{k}),t_n) \dx - \bar\partial_\tau^{\alpha} m_\delta (t_n) +\int_{\partial \Omega} g(t_n) \ds  }{m_\delta(t_n)}\Big],\quad n = 1, \dots , N,
\end{equation*}
where $P_{\mathcal{A}_N}$ denotes a cut-off operation such that for all $v = (v^n)_{n=1}^N$
\begin{equation*}
P_{\mathcal{A}_N} v^n = \min(\max(v^n,0),c_0).
\end{equation*}
The iteration converges to a limit $q_* = (q_*^n)_{n=1}^N$ linearly such that $q_*^n \in [0,c_0]$ for $ n=1,\ldots,N$ and {for any $1<p <\infty$},
\begin{equation}\label{eqn:conv-fully}
  \| [q_{k}^n  -  q^\dag(t_n)]_{n=1}^N \|_{\ell_\lambda^p(\mathbb{R})} \le  (c\lambda^{-\alpha})^{k}\| [q_{0}^n  -  q^\dag(t_n)]_{n=1}^N \|_{\ell_\lambda^p(\mathbb{R})},
\end{equation}
for a sufficiently large $\lambda$. Moreover, the limit $q_*$ satisfies the following error bound
\begin{equation}\label{eqn:inv-err}
  \| [q_*^n  -  q^\dag(t_n)]_{n=1}^N \|_{\ell^p(\mathbb{R})} \le  c (\delta\tau^{-\alpha}+h^2 + \tau^{1/p} |\log \tau|).
\end{equation}
{Here the constant $c$ depends on $\|u_0\|_{H^2(\Omega)}$, $\|g\|_{C^2([0,T];H^\frac12(\partial\Omega))}$, $\|q^\dag\|_{C^1[0,T]}$, Lipschitz constant of $f$, $p$, and $m$, but not on $h$, $\tau$, $n$ and $k$. }
\end{theorem}
\begin{proof}
We define an operator $K_\tau:\mathcal{Q}_N \rightarrow \mathcal{Q}_N$ such that for any $q \in \mathcal{Q}_N$
$$ (K_\tau q)^n = P_{\mathcal{Q}}\Big[\frac{\int_{\Omega}f(u_h^n(q),t_n) \dx
- \bar\partial_\tau^{\alpha} m_\delta (t_n) +\int_{\partial \Omega} g(t_n) \ds  }{m_\delta(t_n)}\Big],\quad n=1,2,\ldots,N.$$
Next, we show that $K_\tau$ is a contraction in $\ell^p_\lambda(\mathbb{R})$ for
a sufficiently large $\lambda$. Let $q_1,q_2 \in \mathcal{Q}_N$. Then $w_h^n = u_h^n(q_1) - u_h^n(q_2)$
satisfies $w_h^n = 0$ and for $n=1,\ldots,N$,
\begin{align}\label{eqn:whn}
   \bar \partial_\tau^\alpha  w_h^n + A_h w_h^n = P_h(f(u_h^{n-1}(q_1),t_n)-f(u_h^{n-1}(q_2),t_n)) + u_h^n(q_2)(q_2^n - q_1^n) -  q_1^n w_h^n.
\end{align}
Then using the fully discrete solution operator $E_{h,\tau}^j$ in \eqref{eqn:FEht-0}, $w_h^n$ can be represented by
\begin{align*}
 w_h^n = \tau \sum_{j=1}^N E_{h,\tau}^{n-j} \big[P_h(f(u_h^{j-1}(q_1),t_j)-f(u_h^{j-1}(q_2),t_j)) + u_h^j(q_2)(q_2^j - q_1^j) -  q_1^n w_h^j\big].
\end{align*}
Using the \textit{a priori} estimate \eqref{eqn:Ehn-stab} and Lipshictz continuity of $f$ in $u$, we obtain
\begin{align*}
 \| w_h^n \|_{L^2\II} &\le  c\tau \sum_{j=1}^n (t_{n+1} - t_j)^{\alpha-1} \big(\| w_h^{j-1}\|_{L^2\II} + \|u_h^j(q_2)\|_{L^2\II}|q_2^j - q_1^j| + |q_1^n| \|w_h^j\|_{L^2\II} \big),
\end{align*}
where $c$ depends on the Lipschitz constant of $f$.
In view of the \textit{a priori} estimate $\|u_h^j(q_2)\|_{L^2\II} \le c$ (cf. Lemma \ref{lem:fully-direct}) and the condition $|q_1^n|\le c_0$ for all $n$, we obtain
\begin{align*}
 \| w_h^n \|_{L^2\II} &\le c \tau \sum_{j=1}^n (t_{n+1} - t_j)^{\alpha-1} \big(\| w_h^{j-1}\|_{L^2\II}+\|w_h^j\|_{L^2\II}  + |q_2^j - q_1^j|\big).
\end{align*}
Multiplying $e^{-\lambda t_n}$ on both sides, taking the ${\ell^p(L^2\II)}$ norm and applying Young's inequality lead to
\begin{align*}
&\quad \| (w_h^n)_{n=1}^N \|_{\ell_\lambda^p(L^2\II)}\\
&\le c  \Big(\tau^{p+1} \sum_{n=1}^N \Big|\sum_{j=1}^n e^{-\lambda (t_n - t_j)}(t_{n+1} - t_j)^{\alpha-1}
e^{-\lambda   t_j }\big(\| w_h^{j-1}\|_{L^2\II}+\|w_h^j\|_{L^2\II}  + |q_2^j - q_1^j|\big)\Big|^p\Big)^\frac1p\\
&\le c \Big(\tau \sum_{j=1}^N  e^{-\lambda t_{j-1}} t_{j}^{\alpha-1}\Big) \big( \| (w_h^n)_{n=1}^N \|_{\ell_\lambda^p(L^2\II)} + \|(q_2^n  - q_1^n)_{n=1}^N \|_{\ell_\lambda^p(\mathbb{R})} \big) \\
&\le c \int_0^T e^{-\lambda s} s^{\alpha-1}\,\ds  \big( \| (w_h^n)_{n=1}^N \|_{\ell_\lambda^p(L^2\II)} + \|(q_2^n  - q_1^n)_{n=1}^N \|_{\ell_\lambda^p(\mathbb{R})} \big).
\end{align*}
From the estimate \eqref{eqn:weight-est}, we obtain
\begin{align*}
 \| (w_h^n)_{n=1}^N \|_{\ell_\lambda^p(L^2\II)} \le c \lambda^{-\alpha} \big( \| (w_h^n)_{n=1}^N \|_{\ell_\lambda^p(L^2\II)} + \|(q_2^n  - q_1^n)_{n=1}^N \|_{\ell_\lambda^p(\mathbb{R})} \big),
\end{align*}
{where $c$ depends on the fractional order $\alpha$ and Lipschitz constant of $f$, but not on $\lambda$, $h$ and $\tau$.}
Now by choosing a sufficiently large $\lambda$ such that $c\lambda^{-\alpha} \le \frac12$, we get
\begin{equation}\label{eqn:weight-est-1}
 \| (w_h^n)_{n=1}^N \|_{\ell_\lambda^p(L^2\II)} \le c \lambda^{-\alpha} \|(q_2^n  - q_1^n)_{n=1}^N \|_{\ell_\lambda^p(\mathbb{R})}.
\end{equation}
Then for $q_1, q_2 \in {\mathcal{Q}}$, we use the stability of the cutoff operation and derive
\begin{equation*}
   |(K_\tau q_1)^n -  (K_\tau q_2)^n| \le \frac{\int_\Omega |f(u_h^n(q_1),t_n) - f(u_h^n(q_2),t_n)| \,\dx}{m_\delta(t_n)}
\le c\| u_h^n(q_1)  -  u_h^n(q_2)\|_{L^2\II}.
\end{equation*}
Taking the $\ell_\lambda^p(\mathbb{R})$ norm on both sides and using the estimate \eqref{eqn:weight-est-1} yield
\begin{equation*}
\begin{split}
   \|((K_\tau q_1)^n -  (K_\tau q_2)^n)_{n=1}^N\|_{\ell_\lambda^p(\mathbb{R})}
&\le  c\| (u_h^n(q_1)  -  u_h^n(q_2))_{n=1}^N\|_{\ell_\lambda^p(L^2\II)} \\
&\le  c\lambda^{-\alpha}\| (q_1^n - q_2^n)_{n=1}^N\|_{\ell_\lambda^p(L^2\II)}.
\end{split}
\end{equation*}
Therefore, for a sufficiently large $\lambda$, $K_\tau$ is a contraction map on the space $\ell^p_\lambda(L^2(\Omega))$. Then Banach fixed point theorem implies that the iteration converges to a limit $(q_*^n)_{n=1}^N \in {\mathcal{Q}_N} $ linearly, i.e.,
\begin{equation*}
\begin{split}
   \|(q_k^n -  q_*^n)_{n=1}^N\|_{\ell_\lambda^p(\mathbb{R})} \le (c\lambda^{-\alpha})^k  \|(q_0^n -  q_*^n)_{n=1}^N\|_{\ell_\lambda^p(\mathbb{R})}.
\end{split}
\end{equation*}
Last, we derive the error estimate for the limit $(q_*^n)_{n=1}^N \in {\mathcal{Q}_N}$. Since it is the fixed point of the operator $K_\tau$, we have
\begin{equation*}
\begin{split}
 & |q_*^n - q^\dag(t_n)| \\
 =& \Big| P_{\mathcal{Q}} \Big[\frac{\int_{\Omega}f(u_h^n(q_*),t_n) \dx - \bar\partial_\tau^{\alpha} m_\delta (t_n) +\int_{\partial \Omega} g(t_n) \ds  }{m_\delta(t_n)}\Big] - q^\dag(t_n)\Big|\\
\le &\Big| \frac{\int_{\Omega}f(u_h^n(q_*),t_n) \dx -   \bar\partial_\tau^{\alpha} m_\delta (t_n) +\int_{\partial \Omega} g(t_n) \ds  }{m_\delta(t_n)}
-\frac{\int_{\Omega}f(u(t_n;q^\dag),t_n) \dx - \partial^{\alpha}_t m (t_n) +\int_{\partial \Omega} g(t_n) \ds  }{m(t_n)}\Big|\\
\le &\Big| \frac{\int_{\Omega}f(u_h^n(q_*),t_n) \dx }{m_\delta(t_n)} -\frac{\int_{\Omega}f(u(t_n;q^\dag),t_n) \dx }{m(t_n)}\Big|
    +\Big| \frac{\bar\partial_\tau^{\alpha} m_\delta (t_n)}{m_\delta(t_n)} -\frac{\partial^{\alpha}_t m (t_n) }{m(t_n)}\Big| \\
   & + \Big| \frac{\int_{\partial \Omega} g(t_n) \ds}{m_\delta(t_n)} -\frac{\int_{\partial \Omega} g(t_n) \ds}{m(t_n)}\Big| = \sum_{j=1}^3 \mathrm{I}_j^n.
\end{split}
\end{equation*}
Condition \eqref{eqn:noise} directly implies
\begin{equation*}
\begin{split}
 \| (\mathrm{I}_3^n)_{n=1}^N \|_{\ell^p(\mathbb{R})} \le c  \| g \|_{C([0,T];L^2(\partial\Omega))} \| (m(t_n) - m_\delta(t_n))_{n=1}^N \|_{\ell^p(\mathbb{R})} \le c \delta.
\end{split}
\end{equation*}
For the term ${\rm I}_2^n$, by the condition \eqref{eqn:noise} and Assumption \ref{ass:stab} (ii), we derive
\begin{equation*}
\begin{split}
{\rm I}_2^n &\le\Big| \frac{m(t_n)\bar\partial_\tau^{\alpha} m_\delta (t_n)-m_\delta(t_n)\partial^{\alpha}_t m (t_n)}{m_\delta(t_n)m(t_n)} \Big| \\
&\le c | (m(t_n)-m_\delta(t_n)) \partial^{\alpha}_t m (t_n) | + c |m(t_n) ( \bar\partial_\tau^{\alpha} m_\delta (t_n) - \partial^{\alpha}_t m (t_n) )|\\
&\le c \delta  \|\partial^{\alpha}_t m \|_{C[0,T]} +  c \| m \|_{C[0,T]}| \bar\partial_\tau^{\alpha} m_\delta (t_n) - \partial^{\alpha}_t m (t_n) | \\
&\le c \delta + c | \bar\partial_\tau^{\alpha} m_\delta (t_n) - \partial^{\alpha}_t m (t_n) |.
\end{split}
\end{equation*}
This and Lemma \ref{lem:num-diff} imply
\begin{equation*}
 \| (\mathrm{I}_2^n)_{n=1}^N \|_{\ell^p(\mathbb{R})} \le c (\delta + \delta \tau^{-\alpha} + \tau^{1/p}|\log \tau|)\le c ( \delta \tau^{-\alpha} + \tau^{1/p}|\log \tau|).
\end{equation*}
Similarly, for the term ${\rm I}_1^n$, the continuous embedding $ C([0,T];H^2\II) \hookrightarrow C([0,T]\times\overline\Omega)$ (for $d=1,2,3$) to $u(q^\dag)$ and the condition \eqref{eqn:noise} yield
\begin{equation*}
\begin{split}
{\rm I}_1^n &\le\Big| \frac{m(t_n)\int_{\Omega}f(u_h^n(q_*),t_n) \dx -m_\delta(t_n)\int_{\Omega}f(u(t_n;q^\dag),t_n) \dx}{m_\delta(t_n)m(t_n)} \Big| \\
&\le c\big| m(t_n)\int_{\Omega}f(u_h^n(q_*),t_n) \dx - m_\delta(t_n)\int_{\Omega}f(u(t_n;q^\dag),t_n) \dx \big| \\
&\le c\Big(\big| (m(t_n) - m_\delta(t_n))\int_{\Omega}f(u(t_n;q^\dag),t_n) \dx \big| + \big|  m(t_n)  \int_{\Omega}f(u_h^n(q_*),t_n) - f(u(t_n;q^\dag),t_n) \dx \big|\Big)\\
&\le c\big( | (m(t_n) - m_\delta(t_n))| + | u_h^n(q_*) - u(t_n;q^\dag)|\big) \\
&\le c \delta + c \| u_h^n(q_*) - u_h^n(q^\dag) \|_{L^2\II} + c \| u_h^n(q^\dag) - u(t_n;q^\dag)\|_{L^2\II}.
\end{split}
\end{equation*}
Then from Lemma \ref{lem:fully-direct}, we deduce
\begin{equation*}
\begin{split}
{\rm I}_1^n  \le c (\delta+h^2+\tau t_n^{\alpha-1}) + c \| u_h^n(q_*) - u_h^n(q^\dag) \|_{L^2\II}.
\end{split}
\end{equation*}
Taking the ${\ell_\lambda^p(\mathbb{R})}$ norm on both sides and applying the estimate \eqref{eqn:weight-est-1} yield
\begin{equation*}
\begin{split}
 \| (\mathrm{I}_1^n)_{n=1}^N \|_{\ell_\lambda^p(\mathbb{R})}
&\le c (\delta+h^2)+ c \tau \|(t_n^{\alpha-1})_{n=1}^N\|_{\ell_\lambda^p(\mathbb{R})} + c \lambda^{-\alpha} \|(q_*^n -  q^\dag(t_n))_{n=1}^N \|_{\ell_\lambda^p(\mathbb{R})}\\
&\le  c (\delta+h^2 + \tau^{\min(1,\alpha+\frac1p)}) + c  \lambda^{-\alpha} \|(q_*^n -  q^\dag(t_n))_{n=1}^N \|_{\ell_\lambda^p(\mathbb{R})}.
\end{split}
\end{equation*}
Combining the preceding estimates yields
\begin{equation*}
\begin{split}
 \| (q_*^n - q^\dag(t_n))_{n=1}^N \|_{\ell_\lambda^p(\mathbb{R})}
 &\le  c (\delta+\delta\tau^{-\alpha} + h^2 + \tau^{1/p}|\log \tau| + \tau^{\min(1,\alpha+\frac1p)}) + c  \lambda^{-\alpha} \|(q_*^n -  q^\dag(t_n))_{n=1}^N \|_{\ell_\lambda^p(\mathbb{R})}\\
&\le c ( \delta\tau^{-\alpha} + h^2 + \tau^{1/p}|\log \tau| ) + c  \lambda^{-\alpha} \|(q_*^n -  q^\dag(t_n))_{n=1}^N \|_{\ell_\lambda^p(\mathbb{R})}.
\end{split}
\end{equation*}
By choosing $\lambda$ sufficiently large and applying the norm equivalence \eqref{eqn:quiv-disc}, we get
\begin{equation*}
 \| (q_*^n - q^\dag(t_n))_{n=1}^N \|_{\ell_\lambda^p(\mathbb{R})}\le  c (\delta\tau^{-\alpha} + h^2 + \tau^{1/p}|\log \tau|).
\end{equation*}
This completes the proof of the theorem.
\end{proof}

\begin{remark}\label{rmk:error-bound}
The error estimate in Theorem \ref{thm:recon} provides useful guidelines to choose the discretization parameters $h$ and $\tau$, by properly balancing the terms, i.e., $\tau^{1/p} \sim \tau^{-\alpha}\delta \sim h^2$, which gives the \textit{a priori} choice $\tau=O(\delta^{{p}/({\alpha p+1})})$ and $h=O(\delta^{{1}/({2(\alpha p+1)})})$. This choice gives the following convergence rate for the approximation $q_*=(q_*^n)_{n=1}^N$:
\begin{equation*}
    \|(q^\dag(t_n) - q_*^n)_{n=1}^{N}\|_{\ell^p(\mathbb{R})}\leq c\delta^{1/(\alpha p +1)} |\log \delta|.
\end{equation*}
That is, the approximation $(q_*^n)_{n=1}^N$ enjoys a H\"{o}lder convergence rate. This agrees with the conditional stability result in Theorem \ref{thm:stab}. Note that the scheme does not incorporate explicit regularization, and the regularizing effect is achieved solely via discretization. The error bound $\tau^{-\alpha}\delta + \tau^{1/p}|\log\tau|$ indicates that either a too large or a too small $\tau$ can lead to large errors in the reconstruction $q_*$. In practical computation, it might also be beneficial to apply filtering to the noisy data $m_\delta$ first in order to ensure stable fractional-order numerical differentiation {\rm(}which also represents the main source of ill-conditioning for the inverse problem{\rm)}.
\end{remark}

\section{Numerical results} \label{sec:numer}
Now we present numerical results to illustrate the theoretical findings. Throughout, the exact data $m$ is generated by solving the direct problem \eqref{eqn:fde} using the Galerkin finite element method in space and backward Euler convolution quadrature in time, cf. Section \ref{sec:fully}, on a finer space-time grid (in order to avoid inverse crime). The noisy data $m_{\delta}$ is generated by
\begin{equation*}
m_{\delta}(t_n) = m(t_n) + \epsilon \zeta(t_n),\quad n=1,\ldots,N,
\end{equation*}
where $\zeta(t_n)$ follows the uniform distribution in $[-1,1]$, $\epsilon\geq0$ denotes the relative noise level, and $\{t_n\}_{n=0}^N$ are the grid points of the fine partition of $[0,T]$. Then, to recover the potential $q^\dag$, we follow the steps in Section \ref{inverse_potential} and design an iterative algorithm based on the fixed point iteration in Theorem \ref{thm:recon}. The iteration is initialized to zero, and the algorithm is run for a maximum of 50 iterations (it is observed to converge within 30 iterations). All the computations are carried out on a personal desktop with MATLAB 2023.

First we present numerical results for a 1D problem with the domain $\Omega = (0,1)$, and $T=0.5$. The problem data $f$, $g$ and $u_0$ are given as follows: $f(u,x,t) = u^2$,  $g(x,t) = 0$, $u_0(x) = 1 + \cos(2\pi x)$. We consider
the following three potentials:
(i) a smooth potential $q_1^\dag$, (ii) a piecewise smooth saw-shaped potential  $q_2^{\dagger}$, and (iii) a discontinuous potential $q_3^\dag$, given respectively by
\begin{equation*}
q_1^{\dagger} = 1 + \cos{5t}, \quad
q_2^{\dagger} = \left\{\begin{aligned}
		& -\tfrac{8}{T}t+2.7,\quad 0\leq t \leq \tfrac{T}{4},\\
		& ~~~~~\tfrac{8}{T}t -1.3,\quad \tfrac{T}{4}\leq t \leq \tfrac{T}{2},\\
  		&-\tfrac{8}{T}t + 6.7 ,\quad \tfrac{T}{2} \leq t \leq \tfrac{3}{4}T, \\
  		& ~~~~~\tfrac{8}{T}t-5.3,\quad  \tfrac{3}{4}T\leq t \leq T,
	\end{aligned}\right.\quad \mbox{and}\quad
q_3^{\dagger} = \left\{\begin{aligned}
		& 2.5,\quad 0\leq t < \tfrac{T}{4},\\
		& 1,\quad \tfrac{T}{4}\leq t < \tfrac{T}{2},\\
  		& 2,\quad \tfrac{T}{2} \leq t < \tfrac{3}{4}T, \\
  		& 1.5,\quad  \tfrac{3}{4}T\leq t \leq T.
	\end{aligned}\right.\quad
\end{equation*}
To measure the accuracy of a reconstruction $q_*=(q_*^n)_{n=1}^N$ for the exact one $q^\dag$, we employ the $\ell^2(\mathbb{R})$ error $e(q_*)$, defined by
$e(q_*)= \|(q^{\dagger}(t_n)- q_*^n)_{n=1}^N\|_{\ell^2(\mathbb{R})}$.
In order to observe the desired convergence rate, for a given noise level $\delta$, we choose the discretization parameters $\tau = {O}(\delta^{2/(1+2\alpha)})$ and $h ={O}(\delta^{1/(2(1+2\alpha))})$, following the theoretical analysis, cf. Remark \ref{rmk:error-bound}. Theorem \ref{thm:recon} (with $p=2$) indicates that the iterative scheme produces a sequence $\left\{q_k\right\}$ of approximations converging linearly to a fixed point $q_*$, and the approximation $q_*$ satisfies an error bound ${O}(\delta^{{1}/({1+2\alpha})})$ (up to a log factor). The experiments below are to complement these theoretical predictions.

\begin{figure}[hbt!]
\centering
\setlength{\tabcolsep}{0pt}
\begin{tabular}{ccc}
\includegraphics[width=.33\textwidth]{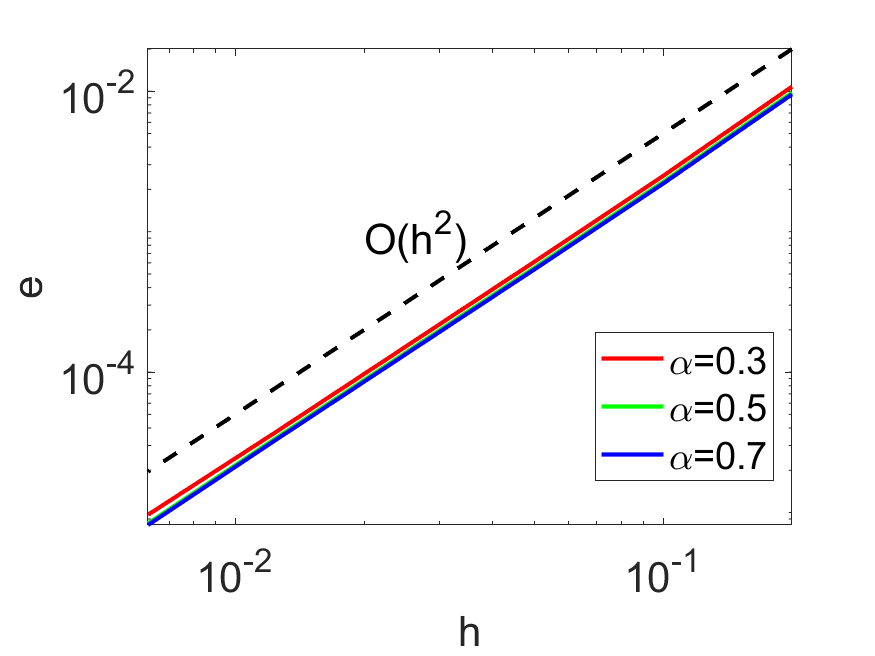} &
\includegraphics[width=.33\textwidth]{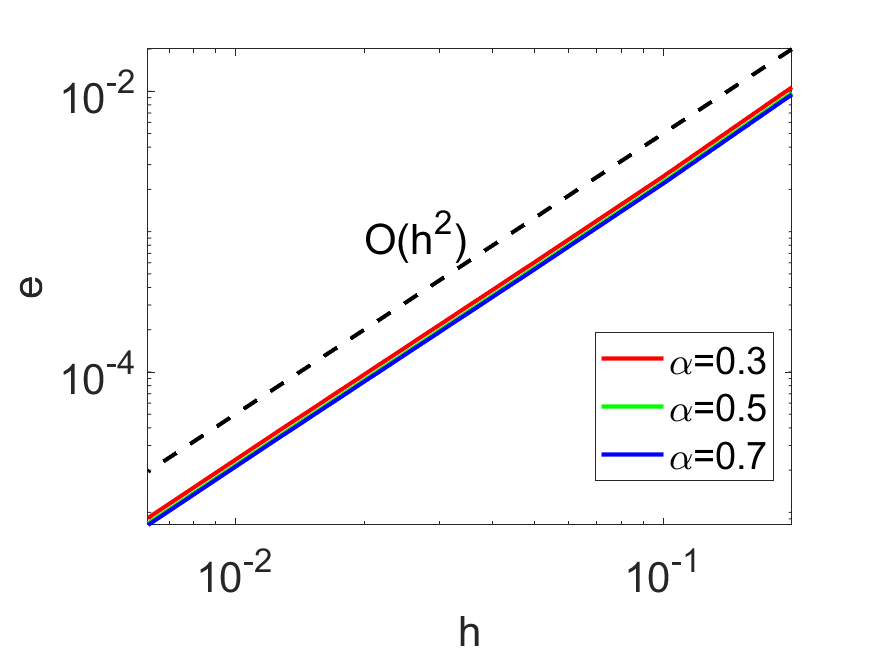} &
\includegraphics[width=.33\textwidth]{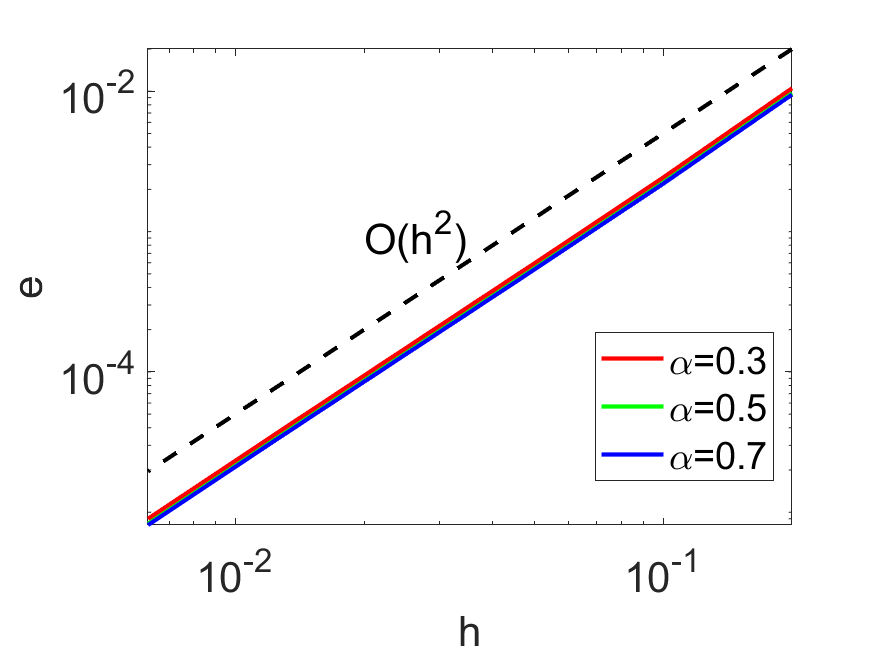} \\
(a) $q_1^\dag$ & (b) $q_2^\dag$ & (c) $q_3^\dag$
\end{tabular}
\caption{The convergence of the approximation with respect to the mesh size $h$, for the 1D case with exact data.\label{fig:conv-h-1d}}
\end{figure}

\begin{figure}[hbt!]
\centering
\setlength{\tabcolsep}{0pt}
\begin{tabular}{ccc}
\includegraphics[width=.33\textwidth]{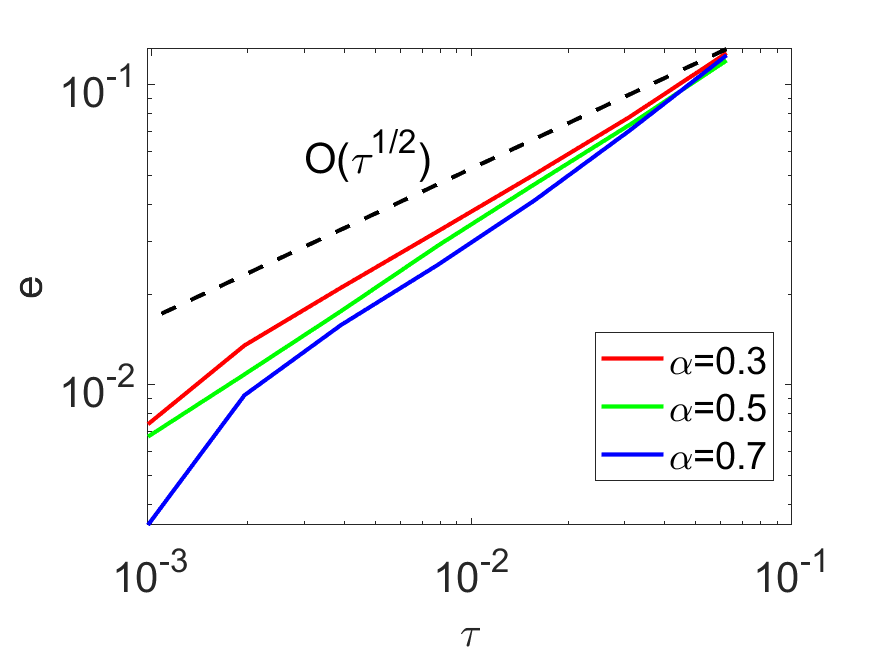} &
\includegraphics[width=.33\textwidth]{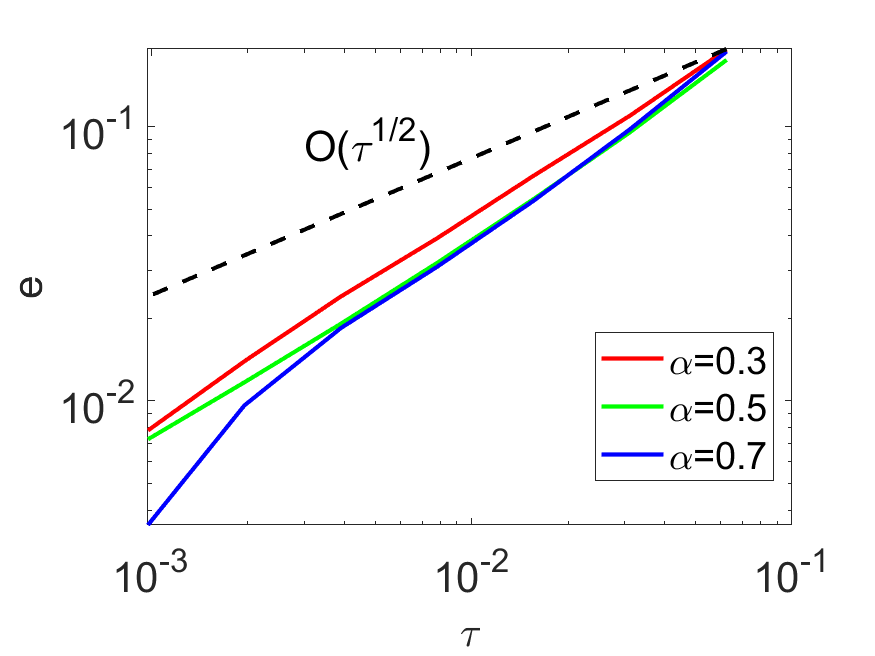} &
\includegraphics[width=.33\textwidth]{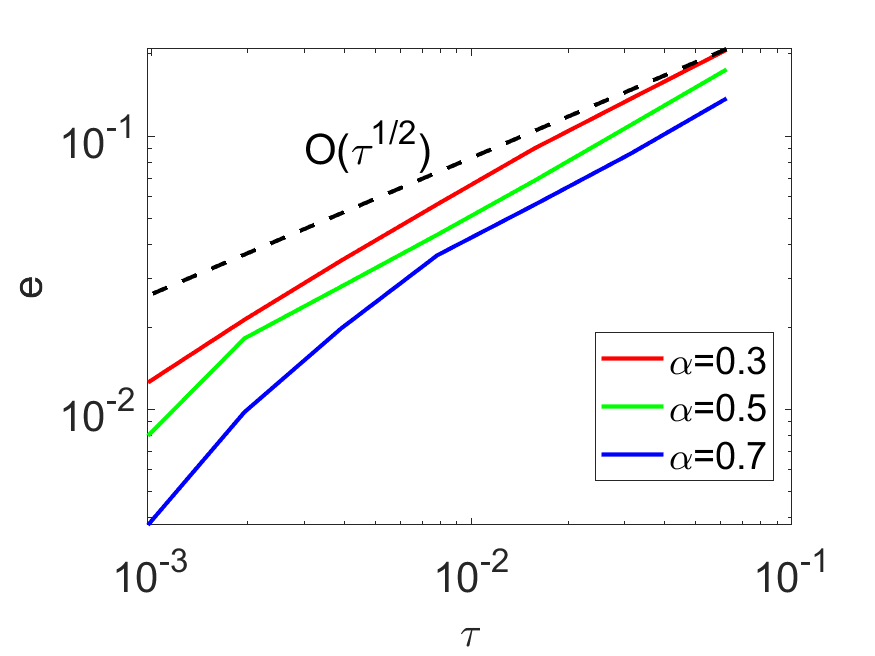} \\
(a) $q_1^\dag$ & (b) $q_2^\dag$ & (c) $q_3^\dag$
\end{tabular}
\caption{The convergence of the approximation with respect to the time step size $\tau$, for the 1D case with exact data.\label{fig:conv-t-1d} }
\end{figure}

First, we numerically test the sharpness of the error estimate in Theorem \ref{thm:recon}:
\begin{equation*}
  \| [q_*^n  -  q^\dag(t_n)]_{n=1}^N \|_{\ell^p(\mathbb{R})} \le  c (\delta\tau^{-\alpha}+h^2 + \tau^{1/p}|\log \tau|).
\end{equation*}
By fixing $\delta = 0$, we examine the impact of the discretization error, which is predicted to be $\mathcal{O}(h^2
+ \tau^{1/p}|\log \tau|)$. To study the spatial convergence, we fix $\tau = T/1000$ and $T=0.5$. The results
in Fig. \ref{fig:conv-h-1d} (for $p=2$) indicate an empirical convergence rate $O(h^2)$ for all potentials
and the fractional order $\alpha$ influences very little the spatial convergence, agreeing well with the
theoretical prediction. Likewise, to investigate the temporal convergence, we fix $h=1\times
10^{-2}$, and $T=0.5$. The results in Fig. \ref{fig:conv-t-1d} exhibit a fairly stable convergence with respect to $\tau$, with an empirical rate slightly
higher than the theoretical one $O(\tau^{1/2}|\log\tau|)$.

\begin{figure}[hbt!]
\centering
\setlength{\tabcolsep}{0pt}
\begin{tabular}{ccc}
\includegraphics[width=.33\textwidth]{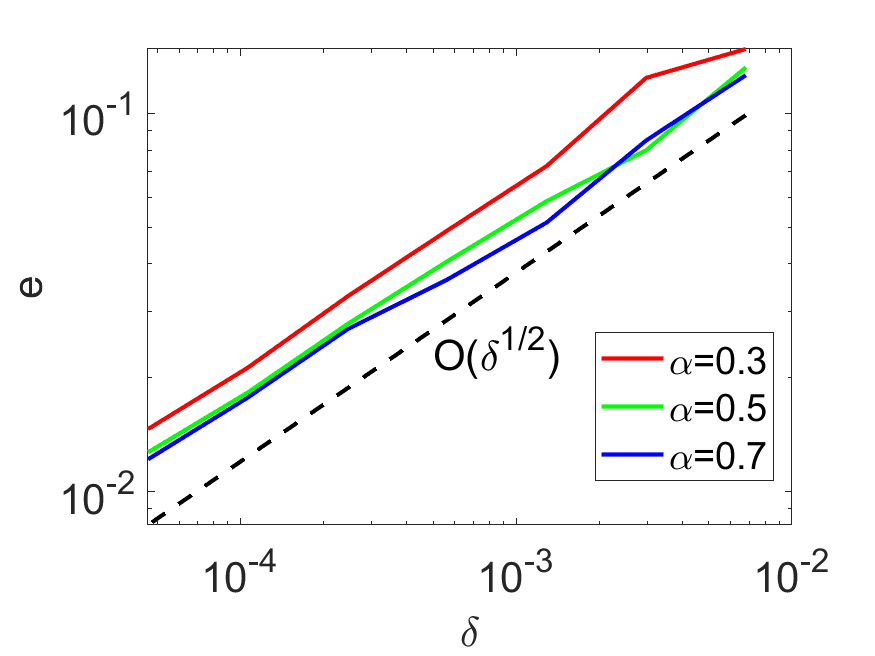} &
\includegraphics[width=.33\textwidth]{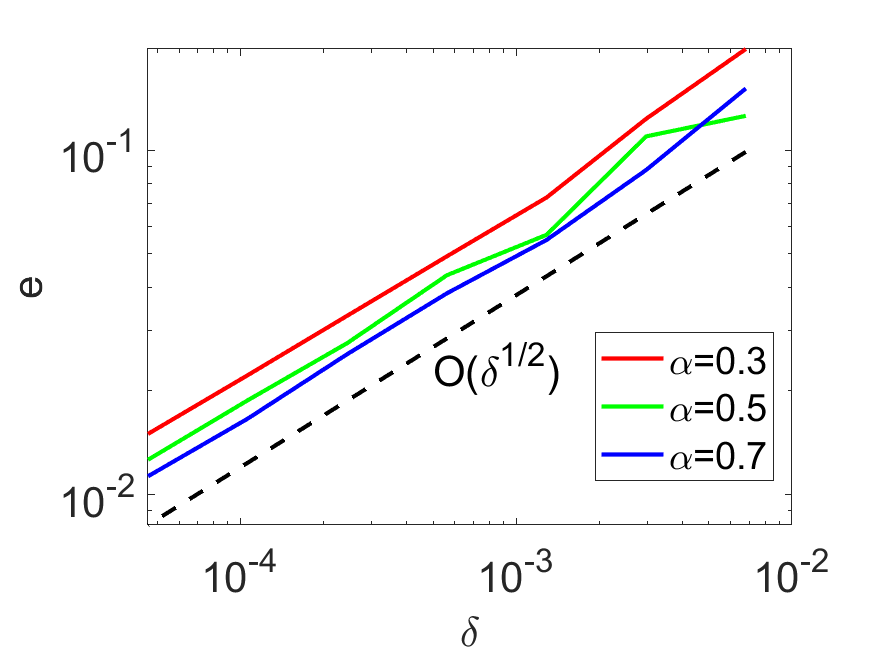} &
\includegraphics[width=.33\textwidth]{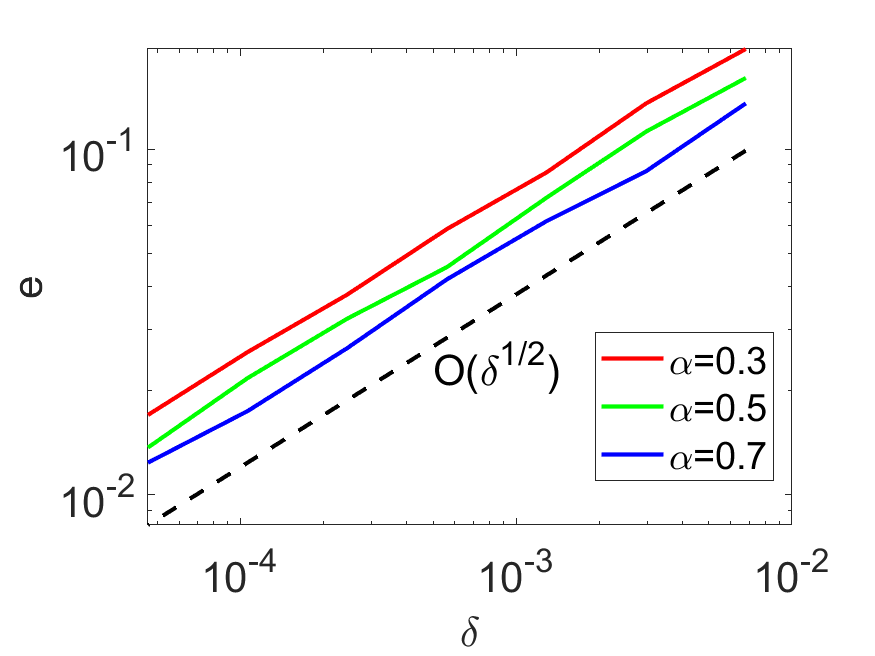} \\
(a) $q_1^\dag$ & (b) $q_2^\dag$ & (c) $q_3^\dag$
\end{tabular}
\caption{The convergence of the approximation with respect to the noise level $\delta$, for the 1D case with noisy data. The approximation is obtained by setting the discretization parameters $h$ and $\tau$ according to Remark \ref{rmk:error-bound}. \label{fig:conv-delta-1d}}
\end{figure}

\begin{figure}
\centering
\setlength{\tabcolsep}{0pt}
\begin{tabular}{ccc}
\includegraphics[width=.33\textwidth]{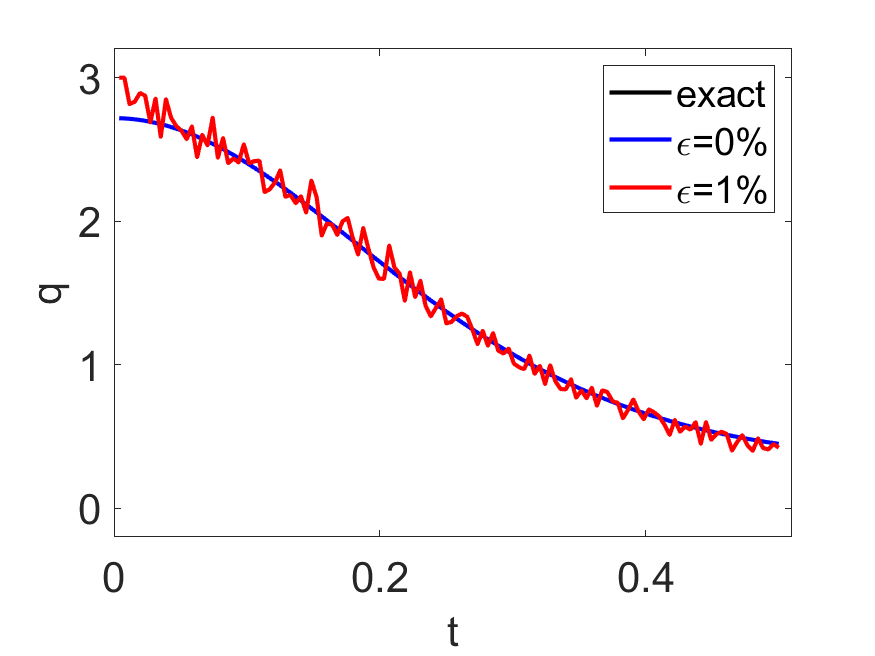} &
\includegraphics[width=.33\textwidth]{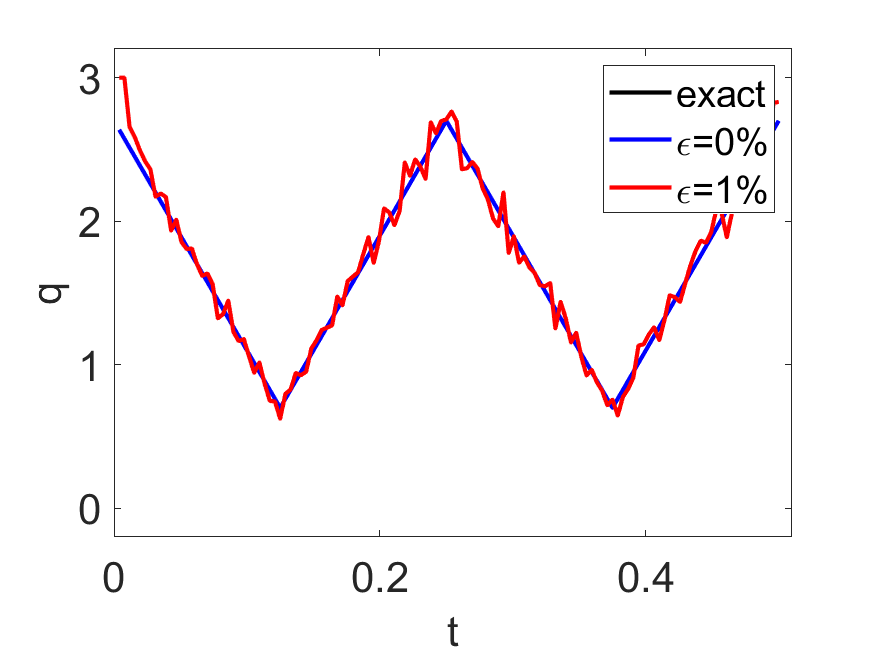} &
\includegraphics[width=.33\textwidth]{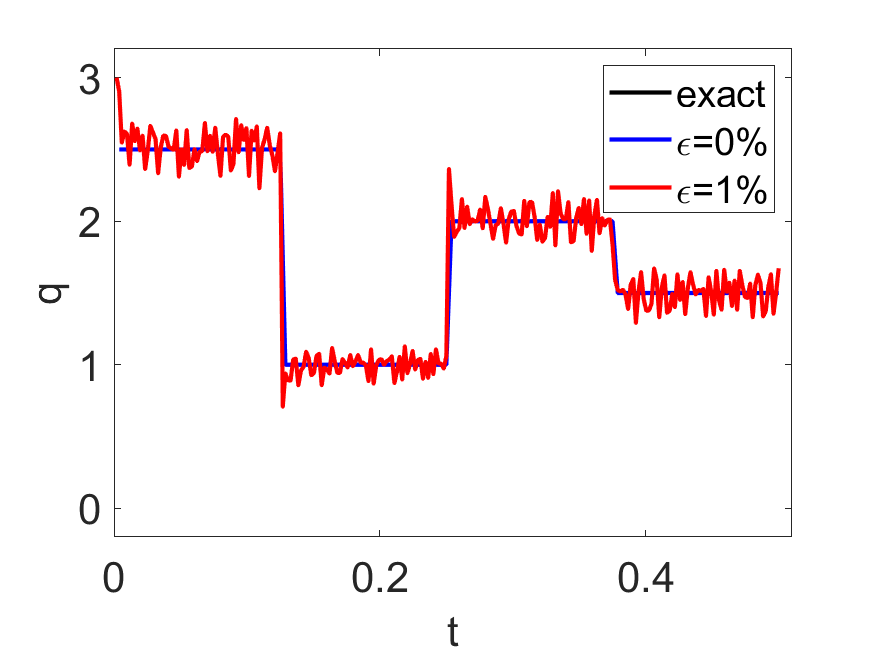}\\
(a) $q_1^\dag $ & (b) $q_2^\dag$ & (c) $q_3^\dag$
\end{tabular}
\caption{The reconstructions by the fixed point iteration for exact and noisy data, with $\alpha=0.5$, and the optimal time step size $\tau$ on the set $\{2^{-i}\}_{i=3}^{11}$.  \label{fig:recon-1d}}
\end{figure}

\begin{figure}[hbt!]
\centering
\setlength{\tabcolsep}{0pt}
\begin{tabular}{cc}
\includegraphics[width=.48\textwidth]{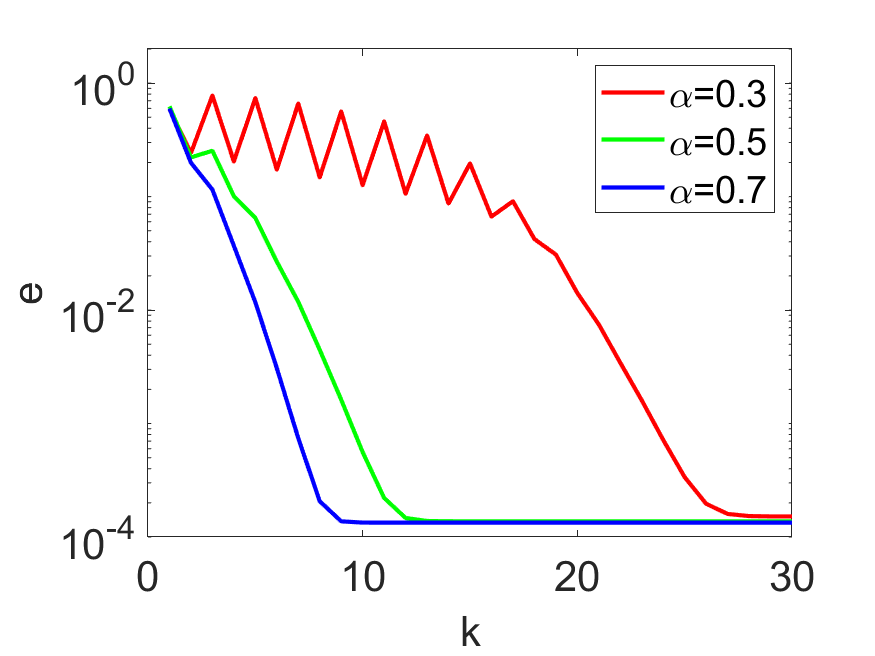} &
\includegraphics[width=.48\textwidth]{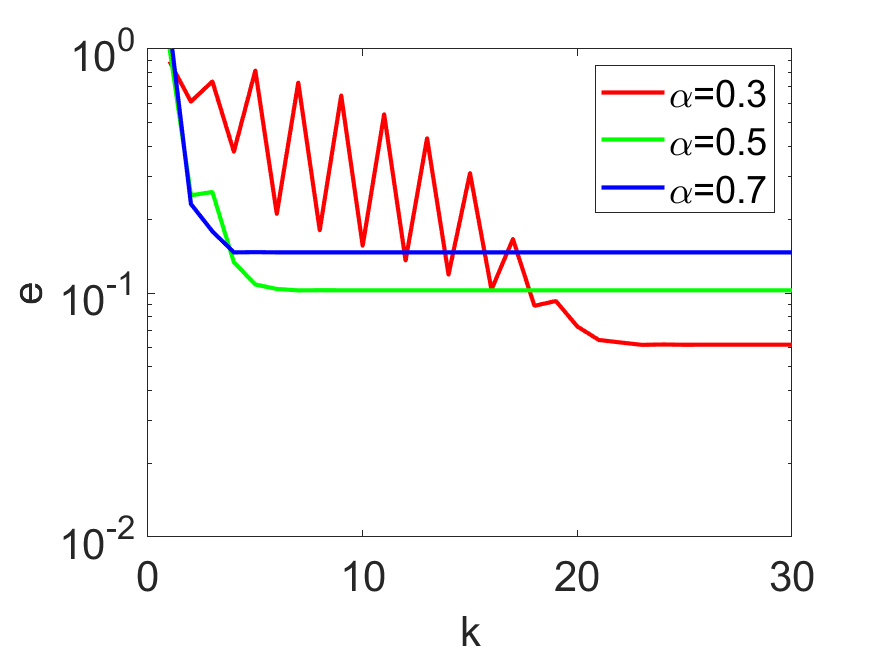} \\
\includegraphics[width=.48\textwidth]{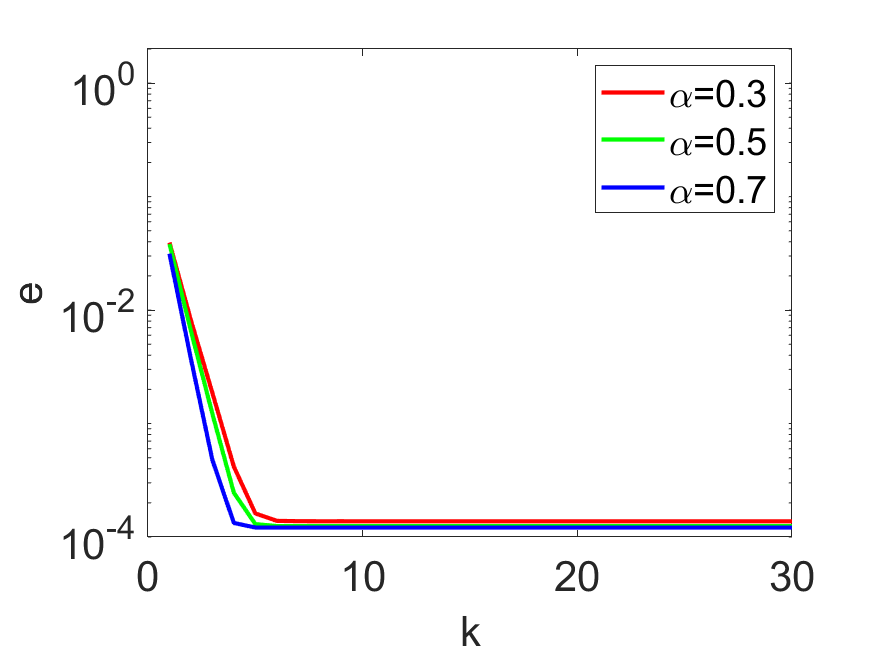} &
\includegraphics[width=.48\textwidth]{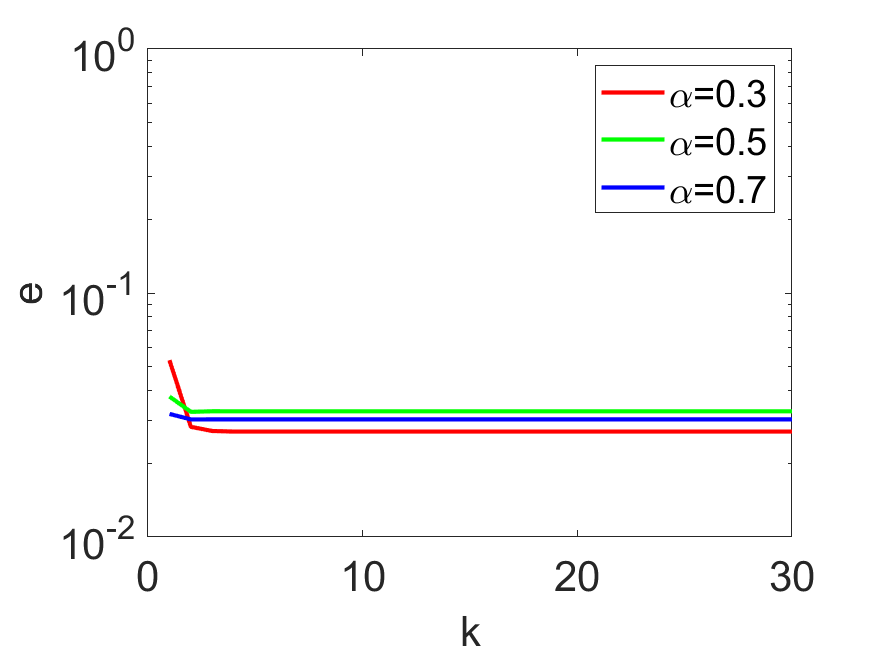} \\
(a) $\epsilon = 0\%$ & (b) $\epsilon = 1\%$
\end{tabular}
\caption{The convergence of the fixed point iteration with respect to the iteration index $k$ for
exact and noisy data. The top row shows the error $e$ in the $\ell^p(\mathbb{R})$, and the bottom row in the $\ell^p_\lambda(\mathbb{R})$ (with $\lambda=25$). The time-discretization level $N$ is optimally chosen to have the best reconstruction accuracy for noisy data.\label{fig:conv-fp}}
\end{figure}

Next, we show the convergence of the reconstructions for noisy observations with
different noise levels $\delta$ in Fig. \ref{fig:conv-delta-1d}, where the results are
obtained by fixing $\tau$ and setting $h =O(\tau^{1/4})$ and $\delta=O(\tau^{\alpha+1/2})$.
The error plots consistently shows an $O(\delta^{1/2})$ empirical convergence rate, which agrees reasonably well with the theoretical prediction $O(\delta^{1/(2\alpha+1)})$
in Theorem \ref{thm:recon}. Interestingly, the observation holds equally well
for the discontinuous potential, which is not covered by the theoretical analysis (since
the admissible set ${\mathcal{Q}}$ only contains continuous potentials). Theoretically, the larger is the fractional order $\alpha$, the more ill-conditioned is \IPP, and hence the
slower is the convergence of the discrete approximation. Nonetheless, numerically, the fractional order $\alpha$ does not appear to influence much the empirical convergence rate, at least for the given parameter choice. Exemplary
reconstructions are shown in Fig. \ref{fig:recon-1d}: for exact data, the numerical recoveries are visually indistinguishable from the exact one, showing the feasibility of the reconstruction, and for noisy data with $\epsilon = 1\%$, the numerical recoveries
exhibit minor oscillations, and the quality of the approximations is largely comparable for all three potentials. The results show that the smoothness of the potential does not play a major in the reconstruction, and the edges / discontinuities can be clearly identified. Note that the small oscillations can be easily removed by filtering \cite{MurioZhan:1998}.

Numerically, the convergence of the
fixed point method in Theorem \ref{thm:recon} does depend very much on the fractional
order $\alpha$: indeed the method converges faster for larger  $\alpha$, which intuitively agrees well with
the convergence analysis in Theorem \ref{sec:recon}, as shown by the factor $O(\lambda^{-\alpha})$
in the estimate \eqref{eqn:conv-fully}. To a certain extent, this observation holds also for noisy data $M_\delta$: the
convergence of the method exhibits a similar behavior but a much smaller number of iterations
are needed in order to reach the convergence, cf. Fig. \ref{fig:conv-fp}. Note that the convergence analysis is performed in the weighted $\ell_\lambda^p(\mathbb{R})$ norm, and it is natural to ask whether it is indeed necessary to use the weighted norm. From Fig. \ref{fig:conv-fp} one observes that the linear convergence holds only for the weighted $\ell^p_\lambda(\mathbb{R})$-norm. The loss of linear convergence in the standard $\ell^p(\mathbb{R})$-norm indeed occurs for small $\alpha$, e.g., $\alpha=0.3$, for which the error plots  clearly exhibit oscillations during the initial stage of the iteration.

\begin{figure}[hbt!]
\centering
\setlength{\tabcolsep}{0pt}
\begin{tabular}{ccc}
\includegraphics[width=.33\textwidth]{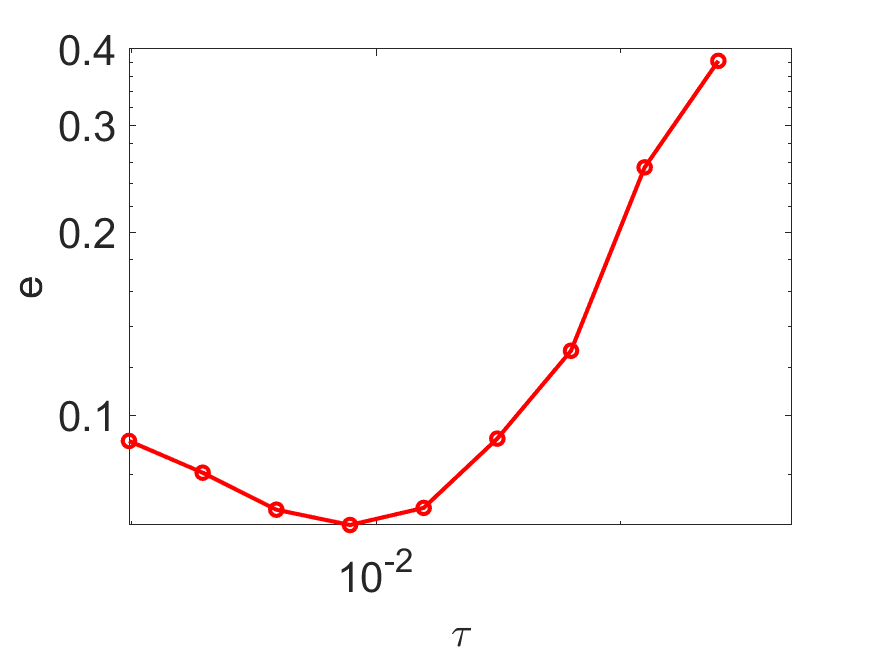} &
\includegraphics[width=.33\textwidth]{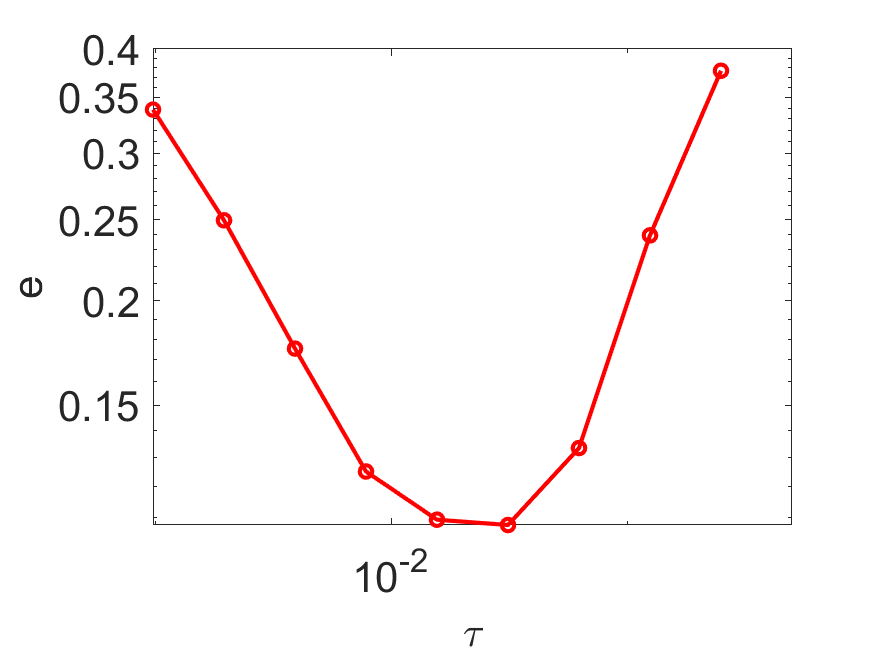} &
\includegraphics[width=.33\textwidth]{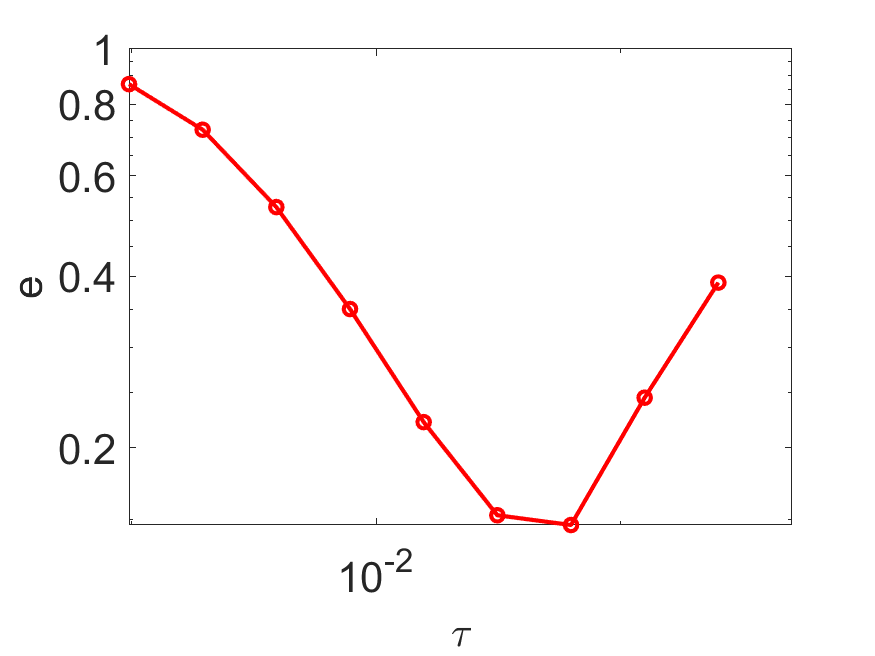}\\
\includegraphics[width=.33\textwidth]{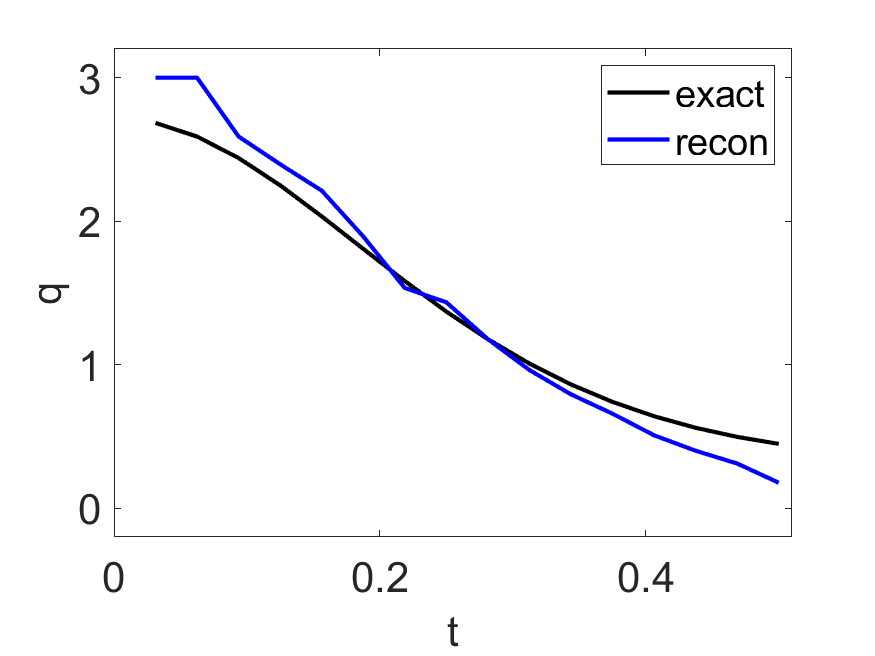} &
\includegraphics[width=.33\textwidth]{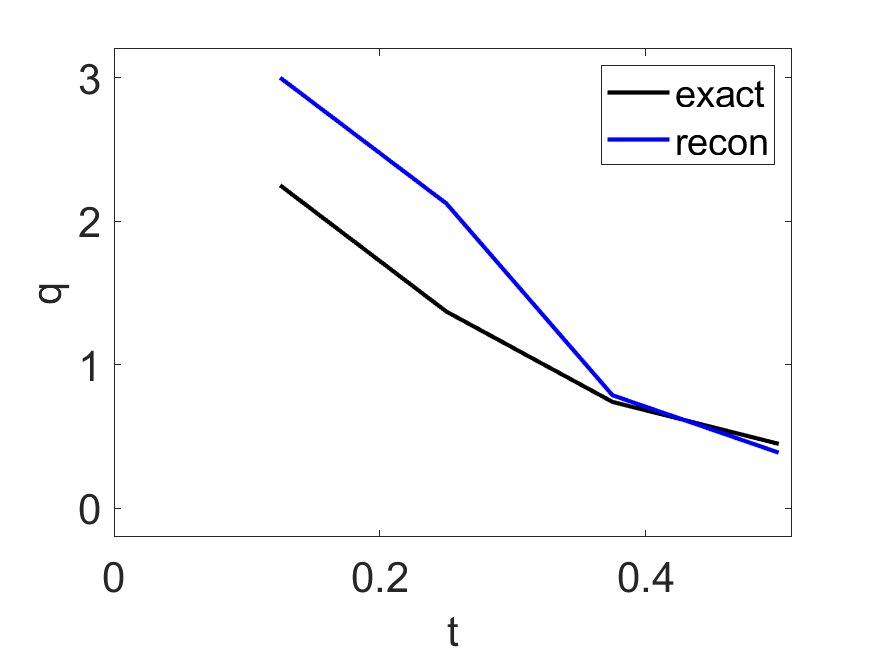} &
\includegraphics[width=.33\textwidth]{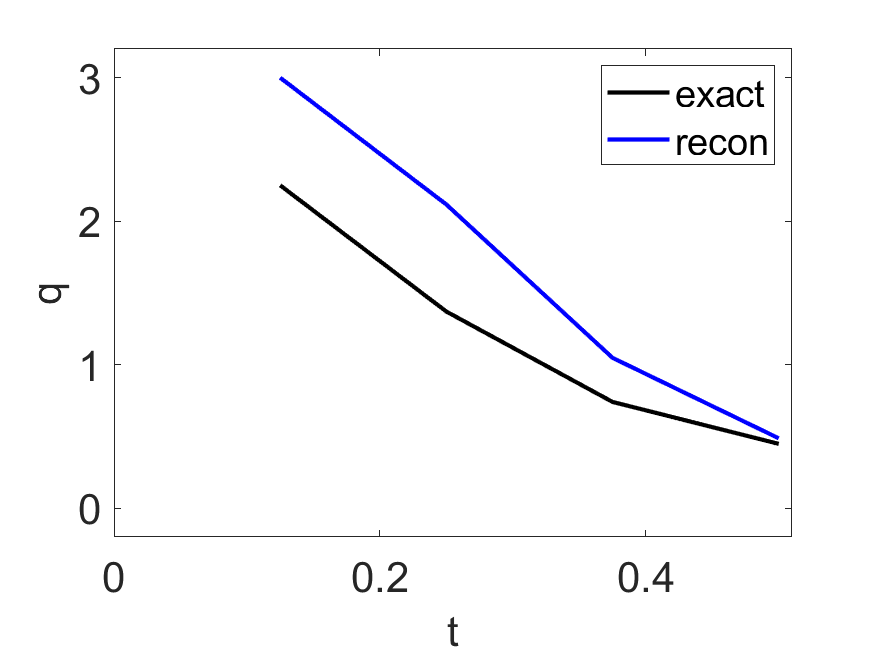} \\
\includegraphics[width=.33\textwidth]{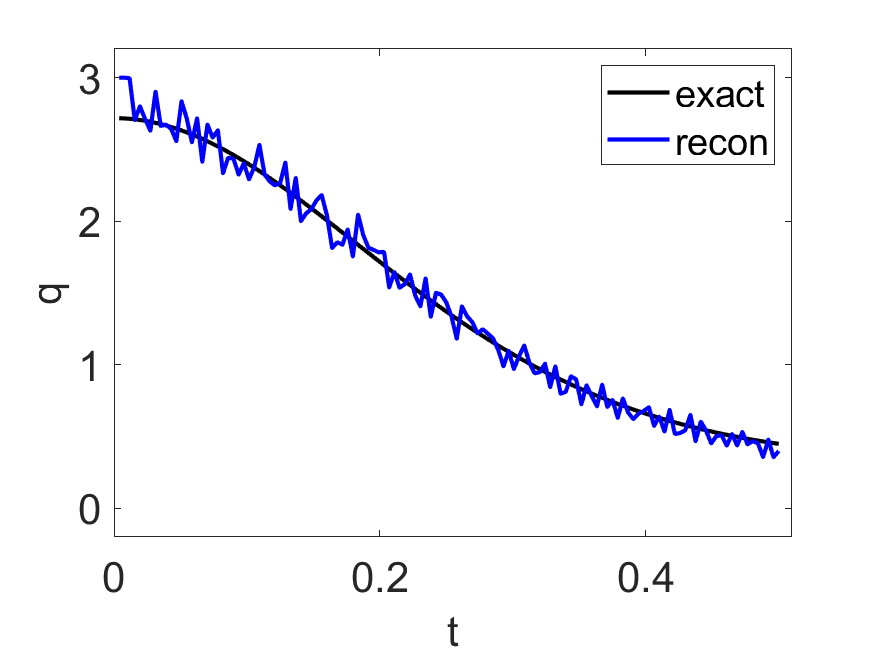} &
\includegraphics[width=.33\textwidth]{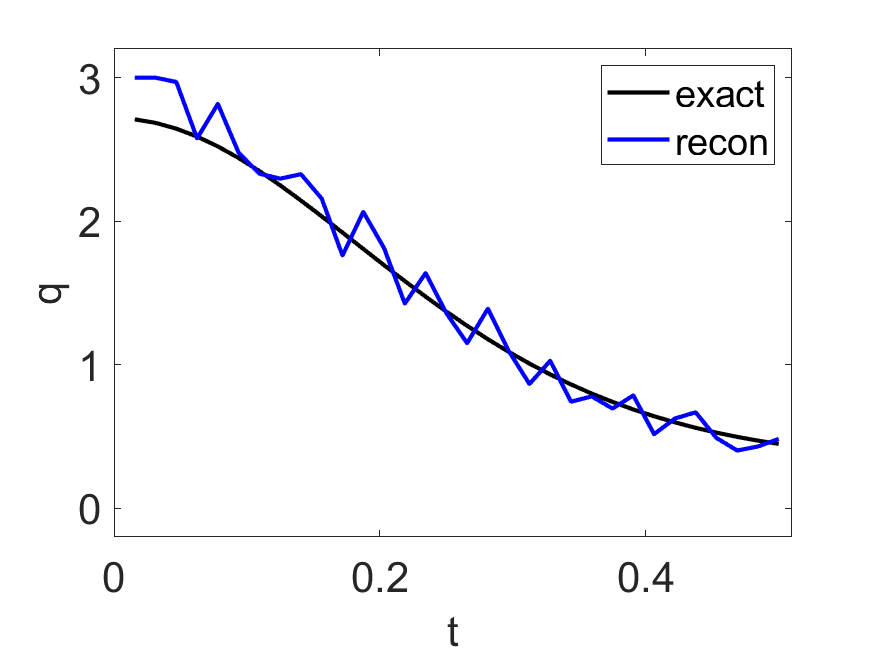} &
\includegraphics[width=.33\textwidth]{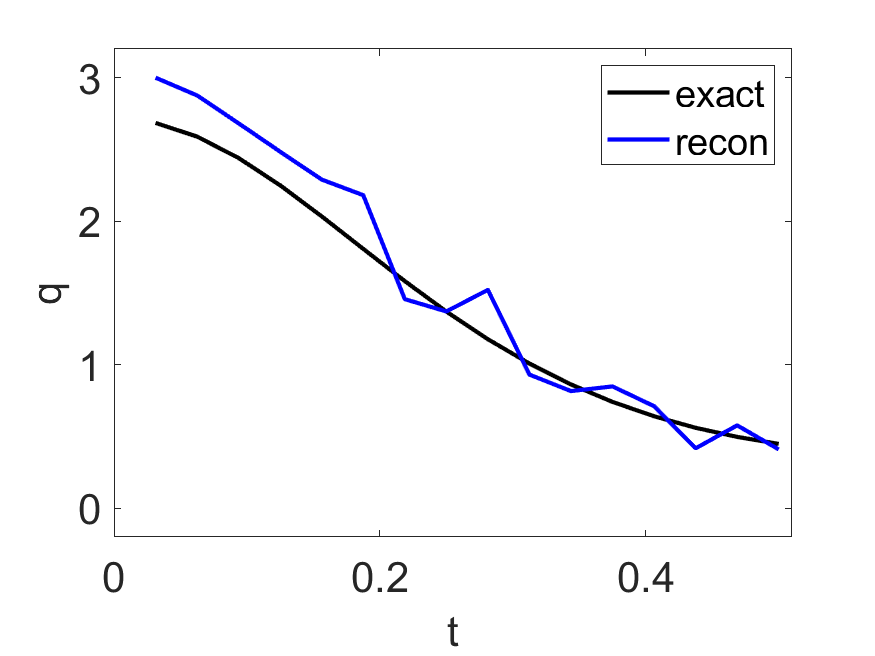} \\
\includegraphics[width=.33\textwidth]{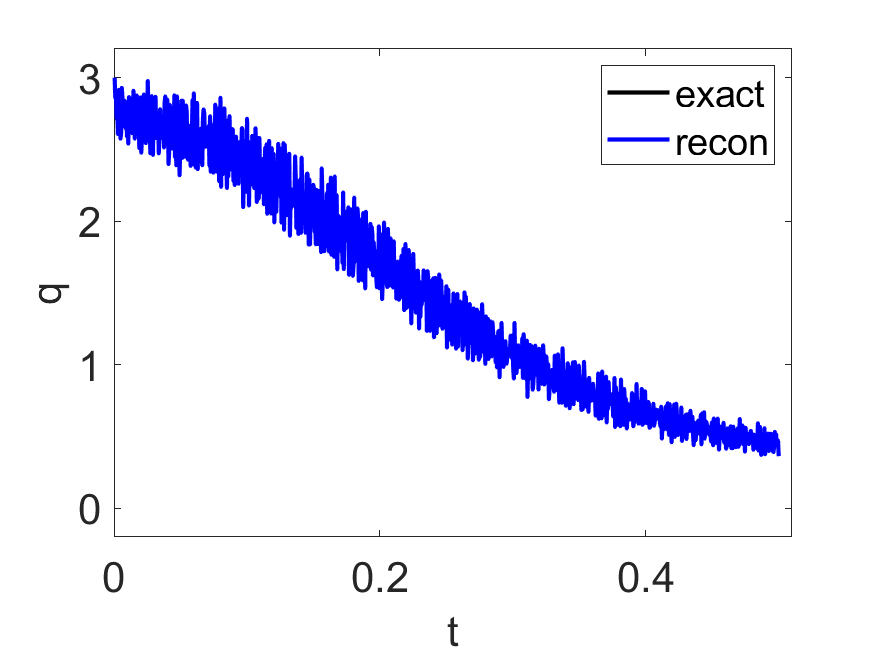} &
\includegraphics[width=.33\textwidth]{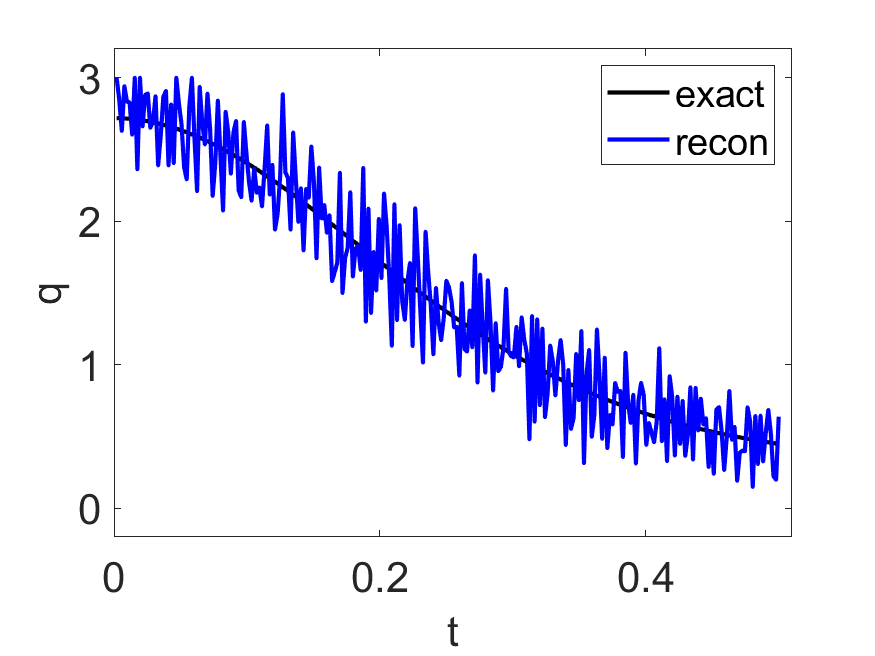} &
\includegraphics[width=.33\textwidth]{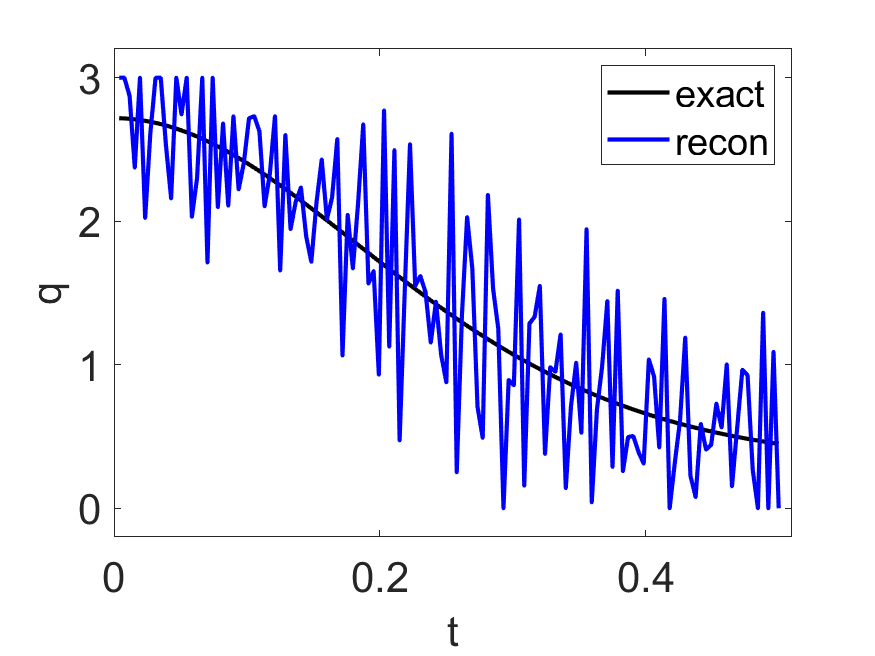} \\
(a) $\alpha=0.3$ & (b) $\alpha=0.5$ & (c) $\alpha=0.7$
\end{tabular}
\caption{The ``semi-convergence'' phenomenon of the fixed point iteration with respect to the time step size, for noisy data with $\epsilon=1\%$,
at three fractional orders, $\alpha=0.3$, $\alpha=0.5$ and $\alpha=0.7$. The top row shows the
$\ell^2(\mathbb{R})$ error of the approximations versus the time step size $\tau$, and the next
three rows show the reconstructions with different time discretization levels. From top to bottom, the total number $N$ of time steps
is $2^{4}$, $2^{7}$ and $2^{10}$ for $\alpha=0.3$; $2^{2}$, $2^{5}$ and $2^{8}$ for
$\alpha=0.5$; and  $2^{2}$, $2^{4}$ and $2^{7}$ for $\alpha=0.7$, where the three discretization levels represent the cases of being over-regularized (too large $\tau$), optimally regularized (optimal $\tau$) and under-regularized (too small $\tau$), respectively. The third row shows the reconstructions with the smallest error. \label{fig:semi-1d} }
\end{figure}

Remark \ref{rmk:error-bound} indicates that the regularizing effect of the fully discrete scheme is achieved by
the time-discretization only, and a suitable choice of the time step size $\tau$ is crucial for the approximation $q_*$ to achieve the best possible accuracy. This behavior is illustrated in Fig. \ref{fig:semi-1d}: in the presence of data noise,
as the step size $\tau$ decreases, the error $e$ first decreases and then increases, indicating
the necessity of choosing an optimal $\tau^*$. When $\tau$ is
optimally chosen, the reconstructions are reasonably accurate (with only mild oscillations),
confirming the conditional stability of \IPP{} in Theorem \ref{thm:stab}. However,  a too large or too small $\tau$ can cause large reconstruction errors, due to the large discretization error (i.e., the factor $\tau^{1/2}|\log \tau|$) or significant amplification of the deleterious
effect of noise (i.e., the factor $\tau^{-\alpha}\delta$). Furthermore, as the fractional order
$\alpha$ increases, IPP becomes increasingly more ill-posed (as indicated by wilder oscillations, resulting
from taking fractional-order differentiation) and thus requires coarser discretizations
(i.e., larger $\tau$); see the top row of Fig. \ref{fig:semi-1d}. In passing, note that the
oscillations can be greatly mitigated by filtering the data a priori or a posterior in order
to obtain visually more appealing approximations.

\begin{figure}[hbt!]
\centering
\setlength{\tabcolsep}{0pt}
\begin{tabular}{ccc}
\includegraphics[width=.33\textwidth]{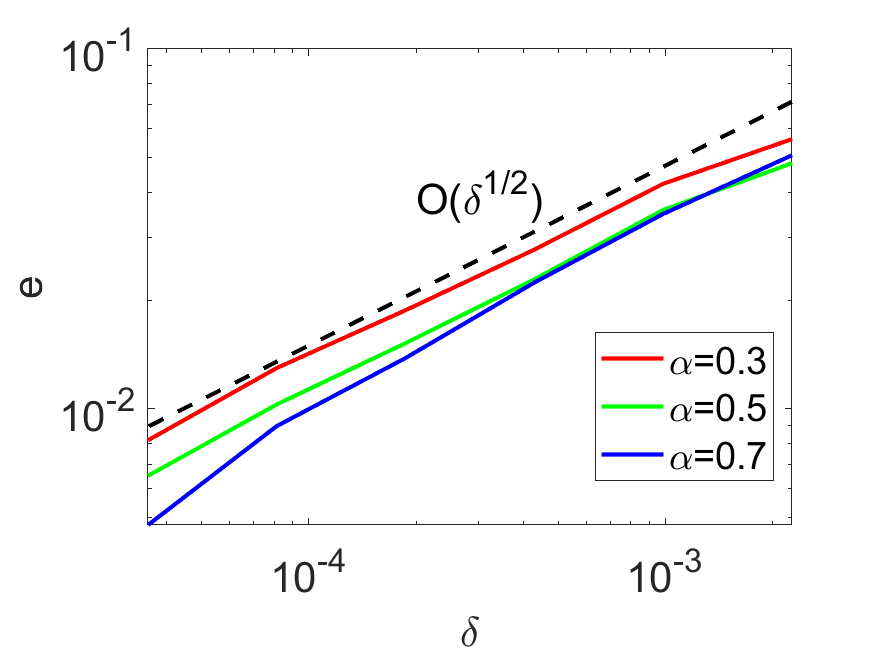} &
\includegraphics[width=.33\textwidth]{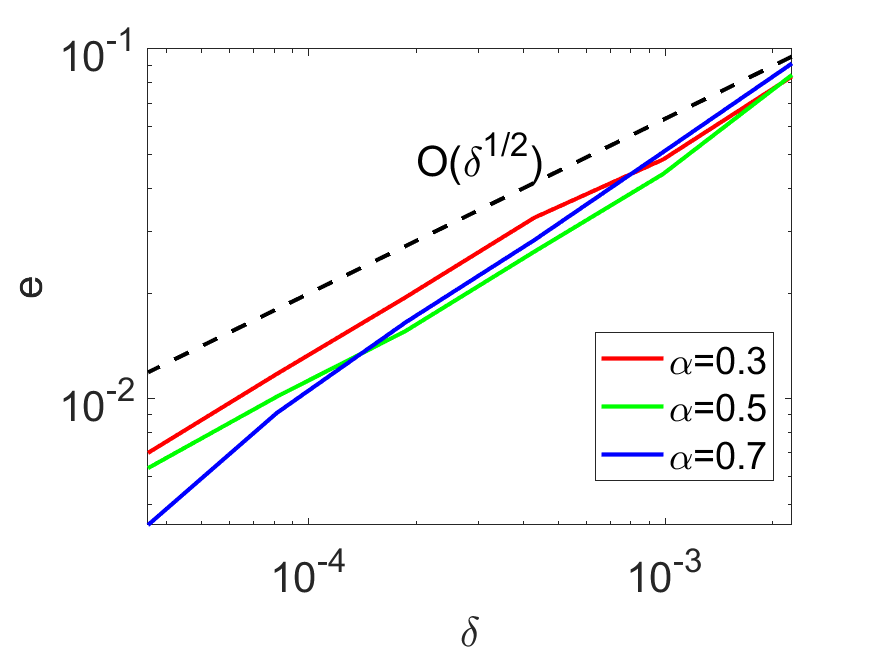} &
\includegraphics[width=.33\textwidth]{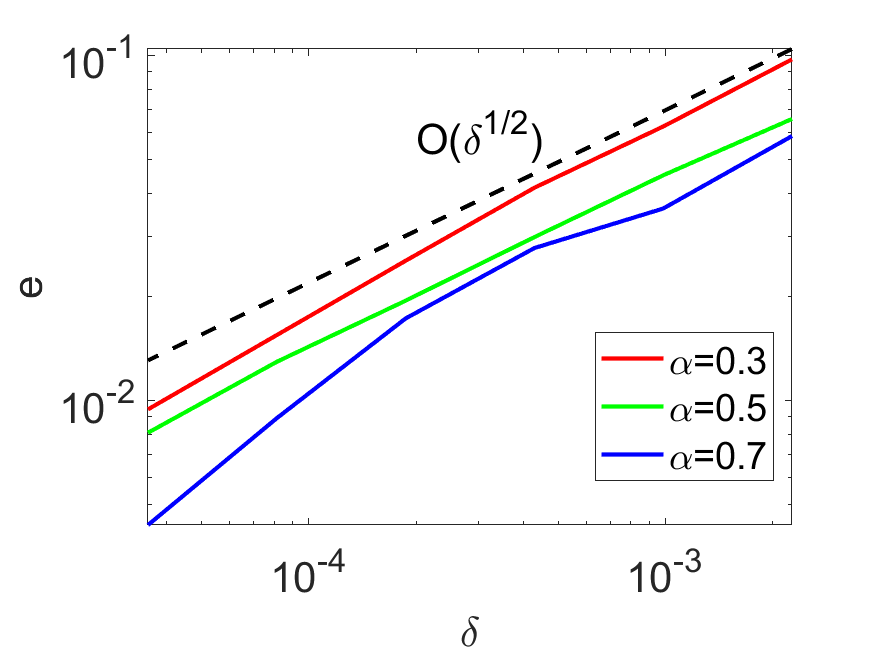} \\
(a) $q_1^\dag$ & (b) $q_2^\dag$ & (c) $q_3^\dag$
\end{tabular}
\caption{The convergence of the approximation with respect to the time step-size $\delta$, for the three potentials in the 2D case with noisy data.\label{fig:conv-delta-2d}}
\end{figure}

Finally, we repeat the experiment with the potentials in the 2D domain $\Omega=(0,1)^2$, with the problem data
$f(u)=u^2$, $g\equiv0$, $u_0(x_1,x_2)=(1+\cos\pi x_1)(1+\cos \pi x_2)$. The numerical results are shown in
Fig. \ref{fig:conv-delta-2d}. The plots show very steady convergence for the discrete approximations, and the empirical convergence rate is close to $O(\tau^{1/2})$, only mildly dependent of the fractional order $\alpha$. These observations agree well with that for the 1D case, which fully agree with the theoretical predictions from Theorem \ref{thm:recon}.

\section{Conclusion}

In this work, we have studied an inverse potential problem of recovering a time-dependent potential $q(t)$ from the spatially averaged measurement over the whole domain. Under minor conditions on the problem data, we have established a conditional Lipschitz stability result, which asserts that it roughly amounts to taking the $\alpha$-order derivative in time of the given data. Furthermore, we developed a fully discrete scheme for approximately recovering the potential $q(t)$, and provided a complete error analysis of the scheme. The obtained error estimate is consistent with the conditional stability estimate. Several numerical experiments were presented to complement the analysis.

{There are several avenues for further research. First, it is of interest to recover the space dependent potential / diffusion coefficient from time-averaged data $\int_0^Tu(x,t)\omega(t){\rm d}t$, for some weight $\omega(t)$. In the standard parabolic case, related inverse problems have been extensively studied (see \cite{PrilepkoKostin:1992} for an early reference). In the fractional case, they are expected to enjoy good stability properties, and error estimates of fully discrete numerical schemes are likely possible. Second, it is also of interest to recover the space (respectively time) dependent potential from the space (respectively time) averaged data $m(t)$. These problems are likely to be much worse behaved, and the analysis would be much more challenging, including numerical analysis of fully discrete schemes. We leave these interesting questions to future works.}

\appendix

\section{Proof of Theorem \ref{thm:sol-reg}}

Now we give a proof of Theorem \ref{thm:sol-reg}, following the fixed point argument of \cite[Theorem 6.17]{Jin:2021}, with minor modifications to treat the nonzero Neumann boundary condition. Using the solution representation
\eqref{eqn:sol} and the Neumann operator $\mathcal{N}$, the solution $u$ satisfies
\begin{equation}\label{eqn:subdiff-semi-re-form}
  u(t) = \mathcal{N}g(t) + F(t)(u_0 - \mathcal{N}g(0)) + \int_0^t E(t-s) [f(u(s),s) - q(s)u(s) - \mathcal{N} \partial_s^\alpha g(s)]\,{\rm d} s.
\end{equation}
The proof is divided into four steps.\\
\noindent{\bf Step 1: Existence and uniqueness.}
We denote by $C([0,T];L^2(\Omega))_\lambda$ the space $C([0,T];L^2(\Omega))$ equipped with the weighted norm
$\|v\|_{\lambda}:=
\max_{0\le t\le T}\|e^{-\lambda t}v(t)\|_{L^2(\Omega)}$, for any $v\in C([0,T];L^2(\Omega)),$
which is equivalent to the standard norm of $C([0,T];L^2(\Omega))$
for any fixed $\lambda>0$. Then we define a map
$M:C([0,T];L^2(\Omega))_\lambda\rightarrow C([0,T];L^2(\Omega))_\lambda$ by
\begin{equation*}
Mv(t)=\mathcal{N}g(t) + F(t)(u_0 - \mathcal{N}g(0)) + \int_0^t E(t-s) [f(u(s),s) - q(s) u(s)- \mathcal{N} \partial_s^\alpha g(s)]\,{\rm d} s.
\end{equation*}
For any $\lambda>0$,
$u\in C([0,T];L^2(\Omega))$ is a solution of problem \eqref{eqn:fde} if and only if it is a fixed point of
the map $M$. It remains to prove that for some $\lambda>0$,
the map $M$ has a unique fixed point.
In fact, the definition of $M$, Lemma \ref{lem:op}(ii), and changing variables $s=t\theta$ yield
\begin{align}
&\quad \|e^{-\lambda t}(Mv_1(t)-Mv_2(t))\|_{L^2(\Omega)} \nonumber\\
&=\bigg\|e^{-\lambda t}\int_0^t E(t-s) [(f(v_1(s),s)-f(v_2(s),s))
+q(s)(v_2(s) - v_1(s))]{\rm d} s\bigg\|_{L^2(\Omega)}\nonumber\\
&\le ce^{-\lambda t}\int_0^t(t-s)^{\alpha-1}\|v_1(s)-v_2(s)\|_{L^2(\Omega)}{\rm d} s\nonumber\\
&\le c\int_0^t(t-s)^{\alpha-1}e^{-\lambda(t-s)}
\max_{s\in[0,T]}\|e^{-\lambda s}(v_1(s)-v_2(s))\|_{L^2(\Omega)}{\rm d} s\nonumber\\
&= c\lambda^{-\alpha}  \bigg(\int_0^1(1-\theta)^{\alpha-1}
(\lambda t)^\alpha e^{-\lambda t(1-\theta)}{\rm d} \theta\bigg) \|v_1-v_2\|_{\lambda}\nonumber \\
&\le c\sup_{\begin{subarray}{c}
\lambda>0,T\geq t>0\\
\theta\in[0,1]
\end{subarray}}
\Big([\lambda t(1-\theta)]^\frac{\alpha}{2}e^{-\lambda t(1-\theta)}\Big)
(\lambda^{-1} t)^\frac{\alpha}{2}
\bigg(\int_0^1(1-\theta)^{\frac\alpha2-1}
{\rm d}\theta\bigg) \|v_1-v_2\|_{\lambda} \nonumber\\
&\le c(\lambda^{-1}T)^\frac{\alpha}{2} \|v_1-v_2\|_{\lambda} ,
\quad\forall\, v_1,v_2\in C([0,T];L^2(\Omega))_\lambda,\label{eqn:subdiff-semi-map-m-contr-1}
\end{align}
{where the constant $c$ depends on the fractional order $\alpha$, $\|q\|_{C[0,T]}$ and Lipschitz constant of $f$.}
By choosing a sufficiently large $\lambda$, we have
 \begin{equation*}
\|e^{-\lambda t}(Mv_1(t)-Mv_2(t))\|_{L^2(\Omega)}
\le \tfrac{1}{2} \|v_1-v_2\|_{\lambda} ,
\quad\forall\, v_1,v_2\in C([0,T];L^2(\Omega))_\lambda .
\end{equation*}
Hence, the map $M$ is contractive on $C([0,T];L^2(\Omega))_\lambda$.
By Banach fixed point theorem, $M$ has a unique fixed point,
which is also the unique solution of problem \eqref{eqn:fde}.

\medskip

\noindent{\bf Step 2: $C^\alpha([0,T];L^2(\Omega))$ regularity.} By the smoothing property of the operator $\mathcal{N}$ and the regularity condition $g\in C^1([0,T];H^\frac{1}{2}(\partial\Omega))$, we have
$\mathcal{N}g \in C^1([0,T];H^2(\Omega))$. Thus, it suffices to analyze the regularity of $w=u-\mathcal{N}g$. Under the given assumption, we have $w_0=u_0-\mathcal{N}g(0)\in D(A)$.
Consider the difference quotient for small $\tau>0$
\begin{align}\label{eqn:subdiff-semi-Holder-u}
\begin{aligned}
\frac{w(t+\tau)-w(t)}{\tau^\alpha}
&=\frac{F(t+\tau)-F(t)}{\tau^\alpha} w_0
+\frac{1}{\tau^\alpha}\int_{t}^{t+\tau} E(s)[f(u(t-s),t-s)-q(t-s)u(t-s)]{\rm d} s \\
  + \int_0^{t} E(t-s) & \frac{f(u(s+\tau),s+\tau)-f(u(s),s) + q(s)u(s) - q(s+\tau)u(s+\tau)}{\tau^\alpha}{\rm d} s =:\sum_{i=1}^3{\mathrm I}_i(t,\tau) .
\end{aligned}
\end{align}
By the identity $F'(t)=-AE(t)$ \cite[Lemma 6.2]{Jin:2021} and Lemma \ref{lem:op} (ii),
\begin{align*}
&\Big\|\frac{F(t+\tau)-F(t)}{\tau^\alpha}w_0\Big\|_{L^2(\Omega)}
 \le \tau^{-\alpha} \int_{t}^{t+\tau}\| F'(s) w_0 \|_{L^2(\Omega)}  \,{\rm d} s\\ \le& c \tau^{-\alpha} \int_t^{t+\tau} s^{\alpha-1}  \,{\rm d} s \|Aw_0\|_{L^2(\Omega)}
\le c \frac{(t+\tau)^{\alpha} - t^\alpha}{\tau^\alpha} \le c.
\end{align*}
By Lemma \ref{lem:op} (ii),
\begin{align*}
\|{\mathrm I}_2(t,\tau)\|_{L^2(\Omega)}
&=\bigg\|\frac{1}{\tau^\alpha}\int_{t}^{t+\tau} E(s)[f(u(t-s),t-s)-q(t-s)u(t-s)]{\rm d} s\bigg\|_{L^2(\Omega)} \\
&\le c\tau^{-\alpha} \int_{t}^{t+\tau}s^{\alpha-1}{\rm d} s= \frac{c}{\alpha} \frac{(t+\tau)^\alpha-t^\alpha}{\tau^\alpha} \le c,
\end{align*}
{where the constant $c$ depends on $c_0$, fractional order $\alpha$ and uniform bound on $f$.}
By the Lipschitz continuity of $f$ and $q$, we have
\begin{align*}
e^{-\lambda t}\|{\mathrm I}_3(t,\tau)\|_{L^2(\Omega)}
&=\bigg\|e^{-\lambda t}\int_0^{t} E(t-s)\frac{f(u(s+\tau),s+\tau)-f(u(s),s)+ q(s)u(s) - q(s+\tau)u(s+\tau)}{\tau^\alpha}{\rm d} s\bigg\|_{L^2(\Omega)} \\
&\le c_1 \int_0^{t}e^{-\lambda (t-s)} (t-s)^{\alpha-1}e^{-\lambda s}\bigg(\bigg\|\frac{w(s+\tau)-w(s)}{\tau^\alpha}\bigg\|_{L^2(\Omega)} + \tau^{1-\alpha}\bigg){\rm d} s  ,
\end{align*}
{where $c$ depends on $\|q\|_{C^1[0,T]}$ and Lipschitz constant of $f$.}
By substituting the bounds on ${\mathrm I}_i(t,\tau)$ into \eqref{eqn:subdiff-semi-Holder-u},
and repeating the argument for \eqref{eqn:subdiff-semi-map-m-contr-1} with the auxiliary function $ W_\tau(t)= e^{-\lambda t}\tau^{-\alpha}\|{w(t+\tau)-w(t)}\|_{L^2(\Omega)}$, we obtain
\begin{equation*}
W_\tau(t)\le c+c_1\int_0^{t} e^{-\lambda (t-s)}(t-s)^{\alpha-1} W_\tau(s){\rm d} s
\le c+c_1(\lambda^{-1}T)^{\frac{\alpha}{2}}\max_{s\in[0,T]}W_\tau(s).
\end{equation*}
By choosing
a large $\lambda$ and taking the maximum in
$t\in[0,T]$, it implies $\max_{t\in[0,T]}W_\tau(t)\le c $, which further yields
\begin{equation*}
  \tau^{-\alpha}\|{w(t+\tau)-w(t)}\|_{X}\le ce^{\lambda t}\le c ,
\end{equation*}
{where the constant $c$ depends also on $T$, but is independent of $\tau$.} Thus, we have proved $\|w\|_{C^\alpha([0,T];X )}
\le c.$\medskip

\noindent{\bf Step 3: $C([0,T];D(A))$ regularity.} By applying the operator $A $ to both sides of
\eqref{eqn:subdiff-semi-re-form} and using the identity $I-F(t)=\int_0^tAE(t-s){\rm d} s$ \cite[Lemmas 6.2 and 6.3]{Jin:2021},
we obtain
\begin{align} \label{eqn:subdiff-semi-Delta-u-reg}
\begin{aligned}
\quad A w(t)-Aw_0&= A  F(t)w_0+\int_0^t A
E(t-s) [f(u(s),s) - q(s)u(s)] {\rm d} s \\
&= \left(AF(t) w_0+ (I-F(t))[f(u(t),t) - q(t)u(t)]\right)\\
&\quad +\int_0^t A  E(t-s)[f(u(s),s)-f(u(t),t) + q(t)u(t) - q(s)u(s) ]{\rm d} s
={\mathrm I}_4(t)+{\mathrm I}_5(t) .
\end{aligned}
\end{align}
By Lemma \ref{lem:op}(i) and the $C^\alpha([0,T];L^2(\Omega))$ regularity from Step 2, we have
\begin{align*}
\|{\mathrm I}_5(t)\|_{L^2(\Omega)}
&=\bigg\|\int_0^t A  E(t-s)[f(u(s),s)-f(u(t),t)+ q(t)u(t) - q(s)u(s)]{\rm d} s\bigg\|_{L^2(\Omega)} \\
&\le c\int_0^t(t-s)^{-1} (\|w(s)-w(t)\|_{L^2(\Omega)} + (t-s)) {\rm d} s\\
&\le c\int_0^t (|t-s|^{\alpha -1} +c) {\rm d} s\le c\max(t,t^\alpha),\quad\forall\, t\in(0,T] ,
\end{align*}
{where the constant $c$ depends on $T$, $\|q\|_{C^1[0,T]}$ and Lipschitz constant of $f$.}
Lemma \ref{lem:op} implies that ${\mathrm I}_5(t)$ is continuous for $t\in(0,T]$,
and the last inequality implies that ${\mathrm I}_5(t)$ is also continuous at $t=0$.
Hence ${\mathrm I}_5 \in C([0,T];L^2(\Omega))$.
Moreover, Lemma \ref{lem:op} gives ${\mathrm I}_4\in C([0,T];L^2(\Omega))$ and
\begin{equation*}
\|{\mathrm I}_4(t)\|_{L^2(\Omega)}\le \|AF(t) w_0+ (I-F(t))[f(u(t),t) - q(t)u(t)]\|_{L^2(\Omega)}\le c.
\end{equation*}
Substituting the bounds on ${\mathrm I}_4(t)$ and ${\mathrm I}_5(t)$ into \eqref{eqn:subdiff-semi-Delta-u-reg} yields
$\|Aw\|_{C([0,T];L^2(\Omega))}\le c$. Then we obtain the estimate $\|w\|_{C([0,T];D(A))}\le c .$
The regularity $u=w+\mathcal{N}g\in C([0,T];D(A) )$ yields
$$\partial_t^\alpha u=Au - qu +f(u,t)\in C([0,T];L^2(\Omega)) .$$

\noindent{\bf Step 4: Estimate of $\|u'(t)\|_{L^2(\Omega)}$.}
Let $q \in C^1[0,T]$. By differentiating
\eqref{eqn:subdiff-semi-re-form} with respect to $t$ and the identity $F'(t)=-AE(t)$ \cite[Lemma 6.2]{Jin:2021}, we obtain (with the subscripts $u$ and $t$ denoting taking partial derivatives in $u$ and $t$, respectively):
\begin{align*}
w'(t) &=F'(t)w_0 + E(t)f(u_0,0)  +\int_0^t E(s) [(f_u(u(s),s) - q(s))u'(s)+(f_s(u,s)-q'(s)u(s))] {\rm d} s \\
&=E(t)(-Aw_0+ f(u_0,0)) +\int_0^t E(t-s) [(f_u(u(s),s) - q(s))u'(s)+(f_s(u,s)-q'(s)u(s))] {\rm d} s .
\end{align*}
By multiplying this equation by $t^{1-\alpha}$, we get
\begin{align*}
t^{1-\alpha}w'(t)
&=t^{1-\alpha}E(t)(-A w_0+f(u_0,0) ) \\
&\quad +\int_0^t t^{1-\alpha}s^{\alpha-1} E(t-s) s^{1-\alpha}
[(f_u(u(s),s)-q(s))  u'(s) + (f_s(u(s),s)-q'(s)u(s))]{\rm d} s ,
\end{align*}
which, together with Lemma \ref{lem:op}(ii), directly implies
\begin{align*}
&e^{-\lambda t}t^{1-\alpha}\|w'(t)\|_{L^2(\Omega)}\\
\le& e^{-\lambda t}t^{1-\alpha}\|E(t)\|_{L^2(\Omega)\to L^2(\Omega)}
\|-Aw_0+ f(u_0,0) \|_{L^2(\Omega)} \\
& +\int_0^t e^{-\lambda (t-s)}t^{1-\alpha}s^{\alpha-1} (t-s)^{\alpha-1} e^{-\lambda s}s^{1-\alpha}(\|f_u(u(s)) - q(s)\|_{L^\infty(\Omega)}
 \|u'(s)\|_{L^2(\Omega)} + \|f_s(u(s),s) - q'(s)u(s) \|_{L^2(\Omega)}){\rm d}s \\
\le &ce^{-\lambda t}\|-A  w_0+ f(u_0,0) \|_{L^2(\Omega) } + c
 + c (\lambda^{-1}T)^{\frac{\alpha}{2}}  \max_{s\in[0,T]}e^{-\lambda s}s^{1-\alpha}\| w'(s)\|_{L^2(\Omega)}  .
\end{align*}
By choosing a sufficiently large $\lambda$ and taking maximum of the left-hand
side with respect to $t\in[0,T]$, we obtain
$\max_{t\in[0,T]}\|e^{-\lambda t} t^{1-\alpha} w'(t)\|_{L^2(\Omega)}
\le c $. This and the bound on $\mathcal{N}g'$ yield the desired bound on $u'$.

\bibliographystyle{abbrv}
\bibliography{frac}

\end{document}